\documentclass[11pt]{article}
\usepackage[utf8]{inputenc}
\usepackage{color}
\usepackage[dvipsnames]{xcolor}
\usepackage{cite}
\usepackage[hidelinks]{hyperref}
\hypersetup{
    colorlinks=true,
    linkcolor=blue,
    citecolor=violet,
    filecolor=magenta,      
    urlcolor=cyan
}
\usepackage[margin = 1in]{geometry}
\usepackage{amsmath, amsfonts,amsthm,amssymb}
\usepackage{graphicx,float,subcaption, soul}

\newcommand{\bmat}[1]{\begin{bmatrix}#1\end{bmatrix}} 

\newcommand{\RR}[1]{\ensuremath{\mathbb{R}^{ #1 }}}
\newcommand{\NN}{\mathbb{N}}
\newcommand{\RRplus}[1]{\ensuremath{\mathbb{R}_+^{ #1 }}}
\newcommand{\defeq}{:=}
\newcommand{\solSymb}{u}
\newcommand{\redsolSymb}{\hat{\solSymb}}
\newcommand{\dotredsolSymb}{\dot{\redsolSymb}}
\newcommand{\snapshotMatSymb}{U}
\newcommand{\funcLibrary}{\Theta}
\newcommand{\DIcoefficients}{\Xi}
\newcommand{\rhsSymb}{f}
\newcommand{\paramSymb}{\mu}
\newcommand{\timeSymb}{t}
\newcommand{\nbig}{N}
\newcommand{\nlibrary}{n_{\ell}}
\newcommand{\nsmall}{n}
\newcommand{\numDI}{n_{\text{DI}}}
\newcommand{\trainingSet}{\mathcal{S}}
\newcommand{\testingSet}{\mathcal{T}}
\newcommand{\setLocalDI}{\trainingSet_{\numDI}}
\newcommand{\indexSet}{\mathcal{I}_{\numDI}}
\newcommand{\localDIindex}{r}
\newcommand{\spaceSymb}{s}
\newcommand{\nat}[1]{\NN(#1)} 

\newcommand{\innat}[1]{\in\nat{#1}}

\newcommand{\paramIndex}{k}
\newcommand{\timeIndex}{j}
\newcommand{\basisSymb}{\Phi}
\newcommand{\basisDummySymb}{\Psi}

\newcommand{\identitySymb}{I}
\newcommand{\leftSingularMatSymb}{V}
\newcommand{\leftSingularVecSymb}{v}
\newcommand{\rightSingularMatSymb}{W}
\newcommand{\singularValueMatSymb}{\Sigma}
\newcommand{\rbfSymb}{\psi}

\newcommand{\rbfcoef}{w}
\newcommand{\rbfcoefk}{\rbfcoef_{k}}

\newcommand{\leftSingularVec}{\boldsymbol{\leftSingularVecSymb}}
\newcommand{\leftSingularMat}{\boldsymbol{\leftSingularMatSymb}}
\newcommand{\rightSingularMat}{\boldsymbol{\rightSingularMatSymb}}
\newcommand{\singularValueMat}{\boldsymbol{\singularValueMatSymb}}
\newcommand{\basis}{\boldsymbol{\basisSymb}}
\newcommand{\basisDummy}{\boldsymbol{\basisDummySymb}}

\newcommand{\snapshotMat}{\boldsymbol{\snapshotMatSymb}}
\newcommand{\reducedSnapshotMat}{\hat{\snapshotMat}}

\newcommand{\sol}{\boldsymbol \solSymb}
\newcommand{\rhs}{\boldsymbol \rhsSymb}
\newcommand{\param}{\boldsymbol{\paramSymb}}
\newcommand{\solinit}{\sol_0}
\newcommand{\finaltime}{\timeSymb_f}
\newcommand{\timeDomain}{[0,\finaltime]}
\newcommand{\paramDomain}{\mathcal D}
\newcommand{\nspacedof}{\nbig_\spaceSymb}
\newcommand{\nreducedspace}{\nsmall_\spaceSymb}
\newcommand{\ntimedof}{{\nbig_\timeSymb}}

\newcommand{\solDummy}{\boldsymbol w}
\newcommand{\timeDummy}{\tau}
\newcommand{\nparam}{\nsmall_\mu}
\newcommand{\dt}{\Delta \timeSymb}
\newcommand{\solapprox}{\tilde{\sol}}
\newcommand{\solapproxFunc}{\solapprox}
\newcommand{\redsolapprox}{\hat\sol}
\newcommand{\redsolapproxFunc}{\redsolapprox}
\newcommand{\identity}{\boldsymbol{\identitySymb}}
\newcommand{\encoder}{\mathcal{G}_{\text{en}}}
\newcommand{\decoder}{\mathcal{G}_{\text{de}}}
\newcommand{\relerr}{e}

\newcommand{\timeArg}[1]{\timeSymb_{#1}}
\newcommand{\region}[1]{\mathcal{D}_{\numDI}(#1)}
\newcommand{\dist}[2]{d({#1},{#2})}
\newcommand{\rhsArg}[1]{\rhs_{#1}}
\newcommand{\solArg}[1]{\sol_{#1}}
\newcommand{\solArgTwo}[2]{\sol_{#1}^{#2}}
\newcommand{\solapproxArg}[1]{\solapprox_{#1}}
\newcommand{\solapproxFuncArg}[1]{\solapproxFunc(t^{#1};\param)}
\newcommand{\redsolapproxArg}[1]{\redsolapprox_{#1}}
\newcommand{\reddotsolapproxArg}[1]{\dot{\redsolapprox}_{#1}}
\newcommand{\redsolapproxArgTwo}[2]{\redsolapprox_{#1}^{#2}}
\newcommand{\redsolapproxFuncArg}[1]{\redsolapproxFunc(\timeSymb^{#1};\param)}
\newcommand{\solFuncArg}[1]{\sol(t^{#1};\param)}
\newcommand{\snapshotMatArg}[1]{\snapshotMat_{#1}}
\newcommand{\identityArg}[1]{\identity_{#1}}
\newcommand{\leftSingularVecArg}[1]{\leftSingularVec_{#1}}
\newcommand{\reducedSnapshotMatArg}[1]{\reducedSnapshotMat_{#1}}

\definecolor{Blue}{rgb}{0,0,1}
\definecolor{Red}{rgb}{1,0,0}
\definecolor{Green}{rgb}{0,1,0}
\definecolor{Cyan}{rgb}{0,0.72,0.92}
\definecolor{Amethyst}{rgb}{0.6,0.4,0.8}
\definecolor{Bronze}{rgb}{0.8,0.5,0.2}
\definecolor{Violet}{rgb}{0.54,0.17,0.89}

\title{LaSDI: Parametric Latent Space Dynamics Identification}
\author{William D. Fries\thanks{Applied Mathematics, School of Mathematical Sciences, University of Arizona, Tucson, AZ 85721 (frieswd@math.arizona.edu)}
\and
Xiaolong He\thanks{Department of Structural Engineering, University of California San Diego, La Jolla, CA 92093 (xih251@eng.ucsd.edu)}
\and 
Youngsoo Choi\thanks{Center for Applied Scientific Computing, Lawrence Livermore National Laboratory, Livermore, CA 94550 (choi15@llnl.gov)}
}

\date{}
\begin{document}

\maketitle






\begin{abstract}
Enabling fast and accurate physical simulations with data has become an important area of computational physics to aid in inverse problems, design-optimization, uncertainty quantification, and other various decision-making applications. This paper presents a data-driven framework for \textit{parametric} latent space dynamics identification procedure that enables fast and accurate simulations. The parametric model is achieved by building a set of local latent space model and designing an interaction among them. An individual local latent space dynamics model achieves accurate solution in a trust region. By letting the set of trust region to cover the whole parameter space, our model shows an increase in accuracy with an increase in training data. We introduce two different types of interaction mechanisms, i.e., point-wise and region-based approach. Both linear and nonlinear data compression techniques are used. We illustrate the framework of Latent Space Dynamics Identification (LaSDI) enable a fast and accurate solution process on various partial differential equations, i.e., Burgers' equations, radial advection problem, and nonlinear heat conduction problem, achieving $O(100)$x speed-up and $O(1)\%$ relative error with respect to the corresponding full order models. 
\end{abstract}



\section{Introduction}\label{sec:intro}
Physical simulations have influenced developments in engineering, technology, and science more rapidly than ever before. The widespread application of simulations as digital twins is one recent example. Advances in optimization and control theory have also increased simulation predictability and practicality. Physical simulations can reproduce responses in a variety of systems, from quantum mechanics to astrophysics, with considerable detail and accuracy, and they provide information that often cannot be obtained from experiments or in-situ measurements due to associated high risks and measurement difficulties. However, high-fidelity forward physical simulations are computationally expensive and, thus, make intractable any decision-making applications, such as inverse problems, design optimization, optimal controls, uncertainty quantification, and parameter studies, where many forward simulations are required.  

To compensate for the computational expense issue, researchers develop several surrogate models to accelerate the physical simulations with high accuracy. One particular type is \textit{projection-based reduced order model} (pROM), in which the full state fields are approximated by applying linear or nonlinear compression techniques. A popular linear compression technique include \textit{proper orthogonal decomposition} (POD) \cite{berkooz1993proper}, reduced basis method \cite{patera2007reduced}, and balanced truncation method \cite{safonov1989schur}, while a popular nonlinear compression technique is \textit{auto-encoder} (AE) \cite{ kim2021fast, kim2020efficient, maulik2021reduced, lee2020model}. The linear compression-based pROM is successfully applied to many different problems, such as Lagrangian hydrodynamics \cite{copeland2022reduced, cheung2022local}, Burgers equations \cite{choi2020sns, choi2019space, carlberg2018conservative},  nonlinear heat conduction problem \cite{hoang2021domain}, aero-elastic wing design problem \cite{choi2020gradient}, Navier--Stokes equation \cite{iliescu2014variational}, computational fluid dynamics simulation for aircraft \cite{amsallem2008interpolation}, convection--diffusion equation \cite{kim2021efficient}, Boltzman transport problem \cite{choi2021space}, topology optimization \cite{choi2019accelerating}, the shape optimization of a simplified nozzle inlet model and the design optimization of a chemical reaction \cite{amsallem2015design}, and lattice-type structure design \cite{mcbane2021component}. These projection-based ROMs are further categorized by its intrusiveness. The \textit{intrusive ROMs} plug the reduced solution representation into the underlying physics law, governing equations, and its numerical discretization methods, such as the finite element, finite volume, and finite difference methods. Therefore, they are physics-constrained data-driven approaches and require less data than the methods that only require data to achieve the same level of accuracy. However, one needs to understand the underlying numerical methods of solving the high-fidelity simulation to implement the intrusive ROMs. Furthermore, the intrusive ROMs are only applicable when the source code for a high-fidelity physics solver is available, which is not the case for certain applications. On the other hand, the non-intrusive ROMs do not require the access to the source code of a high-fidelity physics solver. Therefore, we will focus on developing the \textit{non-intrusive ROMs}, which use \textit{only} data to approximate the full state fields. 

Among many non-intrusive methods, various interpolation techniques are used to build a nonlinear map that predicts new outputs for new inputs. The interpolation techniques include, but not limited to, Gaussian processes \cite{tapia2018gaussian, qian2006building}, radial basis functions \cite{daniel2007hydraulic,huang2015hull}, Kriging \cite{han2013improving,han2012hierarchical}, and convolutional neural networks \cite{guo2016convolutional,zhang2018application}. Among them, neural networks have been the most popular framework because of their rich representation capability supported by the universal approximation theorem. Such surrogate models have been applied to various physical simulations, including, but not limited to, particle simulation \cite{paganini2018calogan}, nanophotonic particle design \cite{peurifoy2018nanophotonic}, porous media flow \cite{kadeethum2021framework,kadeethum2021continuous,kadeethum2021nonintrusive, zhu2018bayesian}, storm prediction \cite{kim2015time}, fluid dynamics \cite{kutz2017deep}, hydrology \cite{marccais2017prospective,chan2018machine}, bioinformatics \cite{min2017deep}, high-energy physics \cite{baldi2014searching}, turbulence modeling \cite{wang2017comprehensive,parish2016paradigm,duraisamy2015new,ling2016reynolds}, uncertainty propagation in a stochastic elliptic partial differential equation \cite{tripathy2018deep}, bioreactor with unknown reaction rate \cite{hagge2017solving}, barotropic climate models \cite{vlachas2018data}, and deep Koopman dynamical models \cite{morton2018deep}. These methods have the potential to be extended to experimental video data as done in \cite{Chen2021}, where intrinsic dimension is sought through auto-encoder training. However, these methods lack the \textit{interpretability} due to the black-box nature caused by its complex underlying structure of interpolators, e.g., neural networks. Furthermore, they do not aim to predict the state field solution in general. Therefore, if the target output quantity of interest changes, the model needs to go through the expensive training phase again, which shows the lack of \textit{generalizability}.

\begin{figure}
    \centering
    \includegraphics[width = \linewidth]{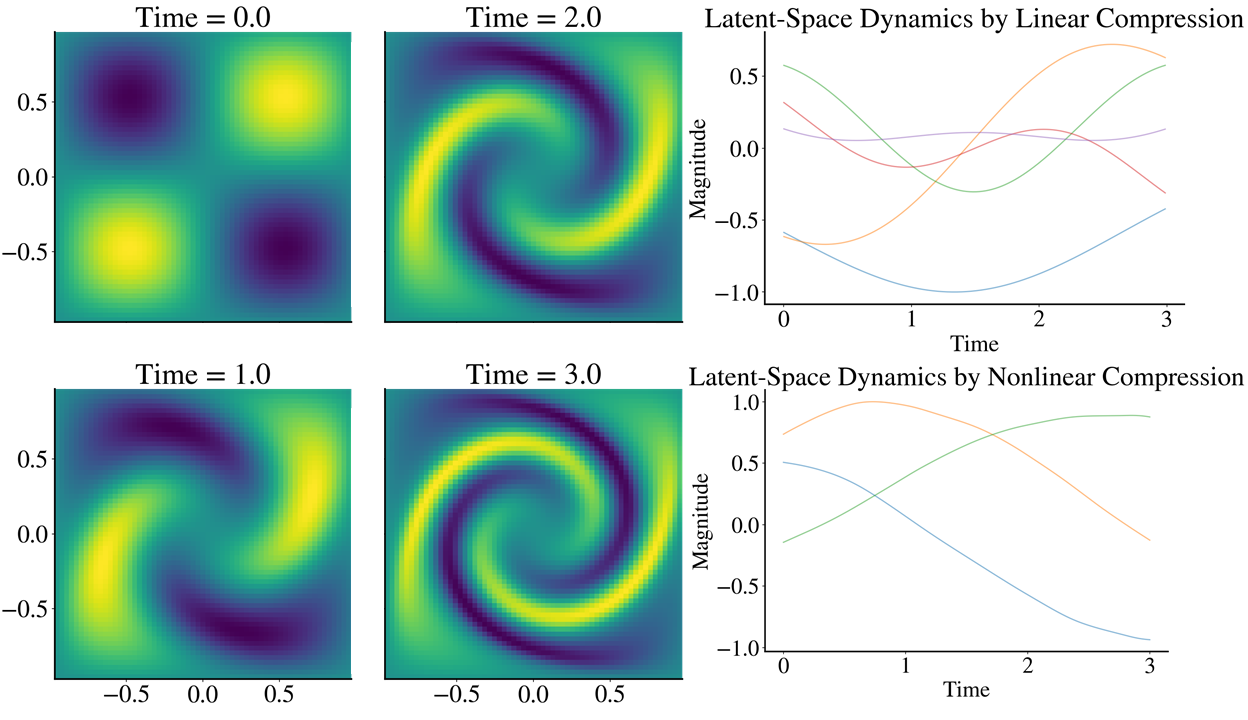}
    \caption{Four snapshots of 2D radial advection simulation are shown, i.e., at $t=0$, $1.0$, $2.0$, and $3.0$. While the full order model involves 9,216 degrees of freedom, the dynamics of the reduced space with five latent variables from a linear compression technique, i.e., proper orthogonal decomposition, and three latent variables from nonlinear compression technique, i.e., auto-encoder, are much simpler than the ones of the full order model.}
    \label{fig:simplicityOfLatentSpaceDynamics}
\end{figure}

To improve the \textit{interpretability} and \textit{generalizability}, we focus on developing a latent space dynamics learning algorithm, where the whole state field data is compressed into a reduced space and the dynamics in the reduced space is learned. Two different types of compression are possible, i.e., linear and nonlinear compression. The popular linear compression can be realized by the singular value decomposition (SVD), justified by the proper orthogonal decomposition (POD) framework. The popular nonlinear compression can be accomplished by the auto-encoder (AE), where the encoder and decoder are designed with neural networks. After the compression, the data size is reduced tremendously. Moreover, the dynamics of the data within the reduced space is often much simpler than the dynamics of the full space. For example, Figure~\ref{fig:simplicityOfLatentSpaceDynamics} shows the reduced space dynamics corresponding to 2D radial advection problems. This motivates our current work to identify the simplified and reduced but latent dynamics with a system identification regression technique. 

There are many existing latent space learning algorithms. DeepFluids \cite{kim2019deep} uses the auto-encoder for nonlinear compression and applies the latent space time integrator to propagate the solution in the reduced space. The authors in \cite{kadeethum2021nonintrusive} and \cite{fresca2021comprehensive} use both linear and nonlinear compressions and apply some interpolation techniques, such as artificial neural networks (ANNs) and radial basis function (RBF) interpolations within the reduced space to predict the solution for new parameter value. Xie, et al., in \cite{xie2019non} uses POD as linear compression technique and apply linear multi-step neural network to predict and propagate the latent space dynamical solutions. Hoang, et al., in \cite{hoang2021projection} compresses both space and time solution space with POD, which was first introduced in \cite{choi2019space}. Then they use several interpolation techniques, such as polynomial regression, $k$-nearest neighbor regression, random forest regression, and neural network regression, within the reduced space formed by the space--time POD. However, all the aforementioned methods use a complex form as a latent space dynamics model, whose structure is not well understood and lacking the interpretability. On the other hand, Champion, et al., in \cite{champion2019data} use an auto-encoder for the compression and the sparse identification of nonlinear dynamics (SINDy) to identify \textit{ordinary differential equations} (ODEs) for the latent space dynamics, which is simpler than neural networks, improving the interpretability. However, they fail to generalize the method to achieve a {\it robust parametric model}. It is partly because one system of ODEs is not enough to cover all the variations within the latent space due to parameter change. 

Therefore, we propose to identify latent space dynamics with a set of local system of ODEs that is tailored to improve the accuracy on a local area of the parameter space. This implies that each local model has a trust region, which covers a sub-region of the whole parameter space. If the set of the trust region covers the whole parameter space, then our model is guaranteed to achieve a certain accuracy level everywhere, achieving truly parametric model based on latent space learning algorithm. At the same time, the framework enables a faster calculation than the corresponding high-fidelity model. We call this framework, LaSDI, which stands for Latent Space Dynamics Identification.

\begin{figure}
    \centering
    \includegraphics[width = \linewidth]{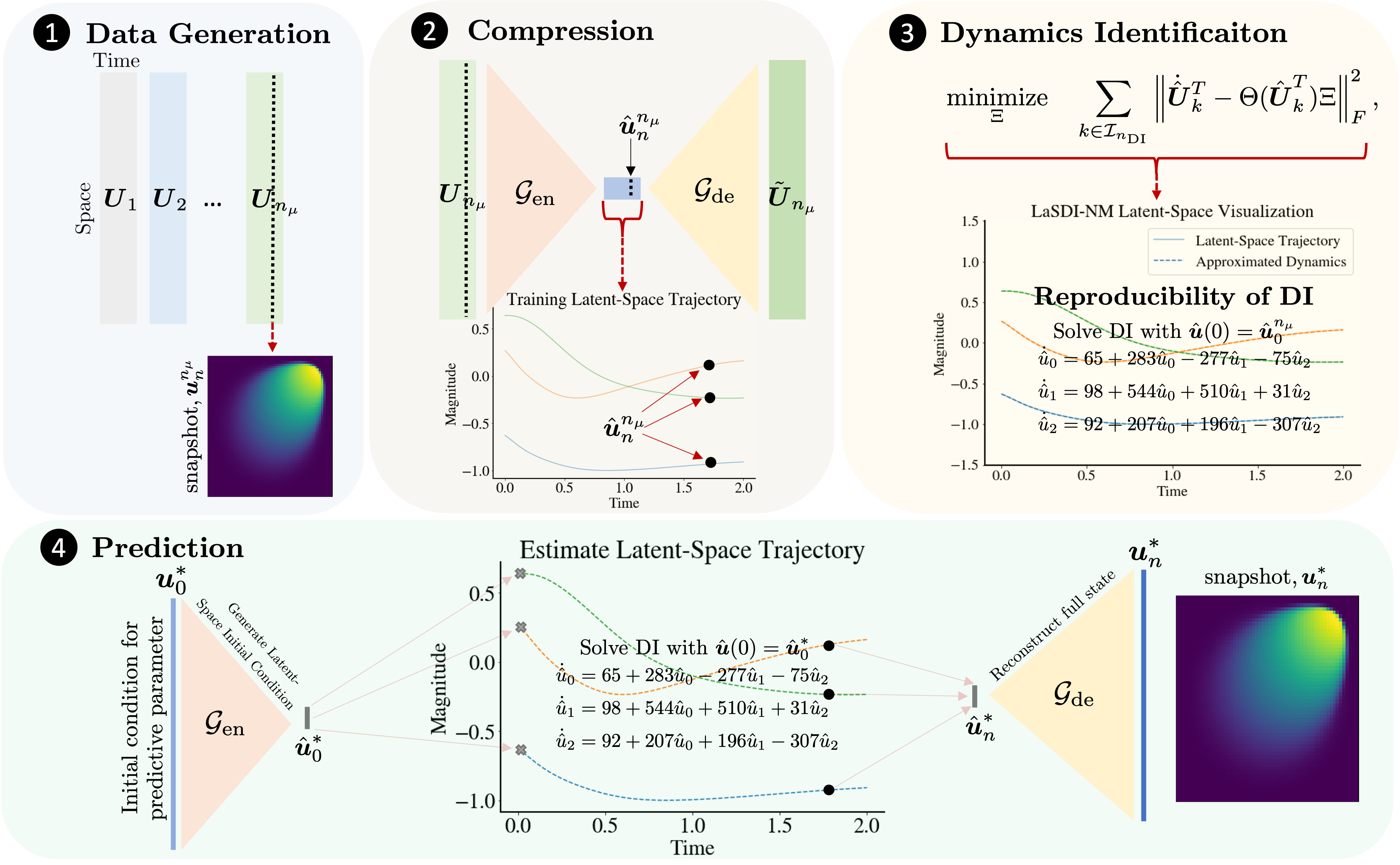}
    \caption{Schematic of LaSDI algorithm applied to 2D Burgers simulations, consisting of four steps: (1) data generation, (2) compression, (3) dynamics identification, and (4) prediction. The full order model simulation data of 2D Burgers problem are collected in Data Generation step (see Section~\ref{sec:data-generation}). Either SVD or autoencoder is used to compress the simulation data to form latent space data in Compression step (see Section~\ref{sec:compression}). Governing equations that best fits the latent space data are identified in Dynamics Identification step (see Section~\ref{sec:dynamicsIndentification}). The identified governing equation is solved with a predictive latent space initial condition and the full state is reconstructed in Prediction  step (see Section~\ref{sec:prediction}).}
    \label{fig:procedureLaSDI}
\end{figure}

The procedure of LaSDI is summarized by four distinct steps below and its schematics are depicted in Figure~\ref{fig:procedureLaSDI}:
\begin{enumerate}
    \item \textbf{Data Generation}: Generate full order model (FOM) simulation data, i.e., parametric time dependent solution data
    \item \textbf{Compression}: Apply compression on the simulation data, through either singular value decomposition or auto-encoder to form latent space data.
    \item \textbf{Dynamics Identification}: Identify the governing equations that best matches the latent space data in a least-squares sense.
    \item \textbf{Prediction}: Use the identified governing equation to predict the latent space solution for a new parameter point. In this paper, we assume that the parameter affects only the initial condition. The predicted latent space solution is reconstructed to the full state by decompression.
\end{enumerate}

In fact, the four-step procedure above is prevalent in developing data-driven reduced-order models to accelerate physical simulations. The traditional data-driven projection-based reduced order models mentioned earlier follow the four-step procedure if the dynamics identification model is \textit{derived} from the existing physical governing equations. The dynamic mode decomposition (DMD) approaches \cite{schmid2010dynamic, tu2013dynamic, proctor2016dynamic, schmid2011applications, kutz2016dynamic, schmid2022dynamic, kutz2016multiresolution, schmid2011application, le2017higher, duke2012error, jovanovic2014sparsity, demo2018pydmd, hemati2014dynamic, erichson2019randomized} also follow the four-step procedure, where linear compression via SVD is used for Step 2 and linear ODEs are generally used for Step 3. Additionally, the operator inference (OpInf) approaches \cite{peherstorfer2016data, geelen2022operator, guo2022bayesian, geelen2022localized, mcquarrie2021non, qian2020lift, swischuk2020learning, benner2020operator, jain2021performance,mcquarrie2021data, peherstorfer2020sampling, khodabakhshi2022non, yildiz2021learning}, which utilize linear compression and fit the latent space dynamics data into polynomial models (mostly quadratic models), also follow the four-step procedure above. As a matter of fact, the proposed LaSDI in this paper is similar to both DMD and OpInf approaches. The major two differences include (i) in Step 2, LaSDI allows for nonlinear compression via an auto-encoder, while DMD and OpInf use a linear compression via the SVD; and (ii) in Step 3, DMD and OpInf only considers polynomial terms in the candidate library, while LaSDI allows general nonlinear terms. Therefore, LaSDI can be viewed as a generalization of both DMD and OpInf.

The main contributions of this paper include
\begin{itemize}
    \item A generalized data-driven physical simulation framework, i.e., LaSDI, is introduced.
    \item A novel local latent space learning algorithm that demonstrates robust parametric prediction, is introduced.
    \item Three novel dynamics identification algorithms are introduced, i.e., global, local, and interpolated ones.
    \item Both good accuracy and speed-up are demonstrated with several numerical experiments. 
\end{itemize}

The paper is organized in the following way: Section~\ref{sec:LaSDI} technically describes LaSDI, where each four distinct steps are explained in details. Section~\ref{sec:data-generation} describes how to generate and store the simulation data. Section~\ref{sec:compression} discusses two different types of compression techniques. Section~\ref{sec:dynamicsIndentification} elaborates the procedure of the dynamics identification and subsequently describes three different types. Section~\ref{sec:prediction} elucidates the prediction step. We demonstrate the performance of LaSDIs for four different numerical experiments in Section~\ref{sec:results}, where global, local, and interpolated dynamics identification algorithms are compared and analyzed. Finally, we summarize and discuss the implications, limitations, and potential future directions for LaSDI in Section~\ref{sec:conclusion}.

\section{Dynamical system of equations}\label{sec:dynamicalSystem}
We formally state a parameterized dynamical system, characterized by the following time dependent ordinary differential equations (ODEs):
\begin{equation}\label{eq:dynamics}
    \frac{d\sol}{dt} = \rhs(\sol,t),\quad\quad
  \sol(0;\param) = \solinit(\param),
\end{equation}
where $\timeSymb\in[0,\finaltime]$ denotes time with the final    time $\finaltime\in\RRplus{}$, and $\sol(\timeSymb;\param)$ denotes the time-dependent, parameterized state implicitly defined as the solution to System~\eqref{eq:dynamics} with $\sol:\timeDomain\times \paramDomain\rightarrow \RR{\nspacedof}$.  Further, $\rhs: \RR{\nspacedof}     \times [0,\finaltime] \rightarrow \RR{\nspacedof}$ with $(\solDummy,\timeDummy)\rightarrow\rhs(\solDummy, \timeDummy) $ denotes the scaled velocity of $\sol$, which can be either linear or nonlinear in its first argument.  The initial state is denoted by $\solinit:\paramDomain\rightarrow \RR{\nspacedof}$, and $\param \in \paramDomain$ denotes the parameters with parameter domain $\paramDomain\subseteq\RR{\nparam}$. We assume that the parameter affects only the initial condition. System~\eqref{eq:dynamics} can be considered as semi-discretized version of a system of partial differential equations (PDEs), whose spatial domain is denoted as $\Omega\in\RR{d}$, $d\innat{3}$, where $\nat{N}\defeq\{1,\ldots,N\}$. The spatial discretization can be done through many different numerical methods, such as the finite difference, finite element, and finite volume methods. 

Many different time integrators are available to approximate the time derivative term, $d\sol/d\timeSymb$, e.g., explicit and implicit time integrators. A uniform time discretization is assumed throughout the paper, characterized by time step $\dt\in\RRplus{}$ and time instances $\timeArg{n} = \timeArg{n-1} + \dt$ for $n\innat{\ntimedof}$ with $\timeArg{0} = 0$, $\ntimedof\in\NN$.  To avoid notational clutter, we introduce the following time discretization-related notations: $\solArg{n} \defeq \solFuncArg{n}$, $\solapproxArg{n} \defeq\solapproxFuncArg{n}$, $\redsolapproxArg{n} \defeq \redsolapproxFuncArg{n}$, and $\rhsArg{n} \defeq \rhs(\solFuncArg{n},t^{n}; \param)$.

Although many advanced numerical methods are available to solve System~\eqref{eq:dynamics}, as the problem size increases, i.e., a large $\nspacedof$, and the computational domain $\Omega$ gets geometrically complicated, the overall solution time becomes impractically slow, e.g., taking a day or week for one forward simulation. The proposed method, i.e., LaSDI, accelerates the computationally expensive simulations. 

\section{LaSDI}\label{sec:LaSDI}
This section mathematically describes LaSDI. As shown in Introduction~\ref{sec:intro}, LaSDI consists of four steps, i.e., data generation, compression, identification, and prediction. Each step is described in the subsequent sections.

\subsection{Data generation}\label{sec:data-generation}
LaSDI first generates full order model simulation data by solving a system of dynamical system \eqref{eq:dynamics} by sampling parameter space, $\paramDomain$. We denote the sampling points by $\param_{\paramIndex}\in\trainingSet\subset\paramDomain$, $\paramIndex\innat{\nparam}$ where $\trainingSet$ denotes a training set. We also denote $\solArgTwo{n}{\paramIndex}\in\RR{\nspacedof}$ for the $n$-th time step solution of \eqref{eq:dynamics} with $\param = \param_{\paramIndex}$ and arrange the snapshot matrix $\snapshotMatArg{\paramIndex} = \bmat{\solArgTwo{0}{\paramIndex} & \cdots & \solArgTwo{\ntimedof}{\paramIndex}} \in \RR{\nspacedof \times (\ntimedof+1)}$ for $\param = \param_{\paramIndex}$. Concatenating all the snapshot matrices side by side, the whole snapshot matrix $\snapshotMat\in\RR{\nspacedof \times (\ntimedof+1)\nparam}$ is defined as
\begin{equation}\label{eq:snapshotMat}
    \snapshotMat = \bmat{\snapshotMatArg{1} & \cdots & \snapshotMatArg{\nparam}}.
\end{equation}

\subsection{Compression}\label{sec:compression}
\begin{figure}
\centering
    \includegraphics[width = .6\textwidth]{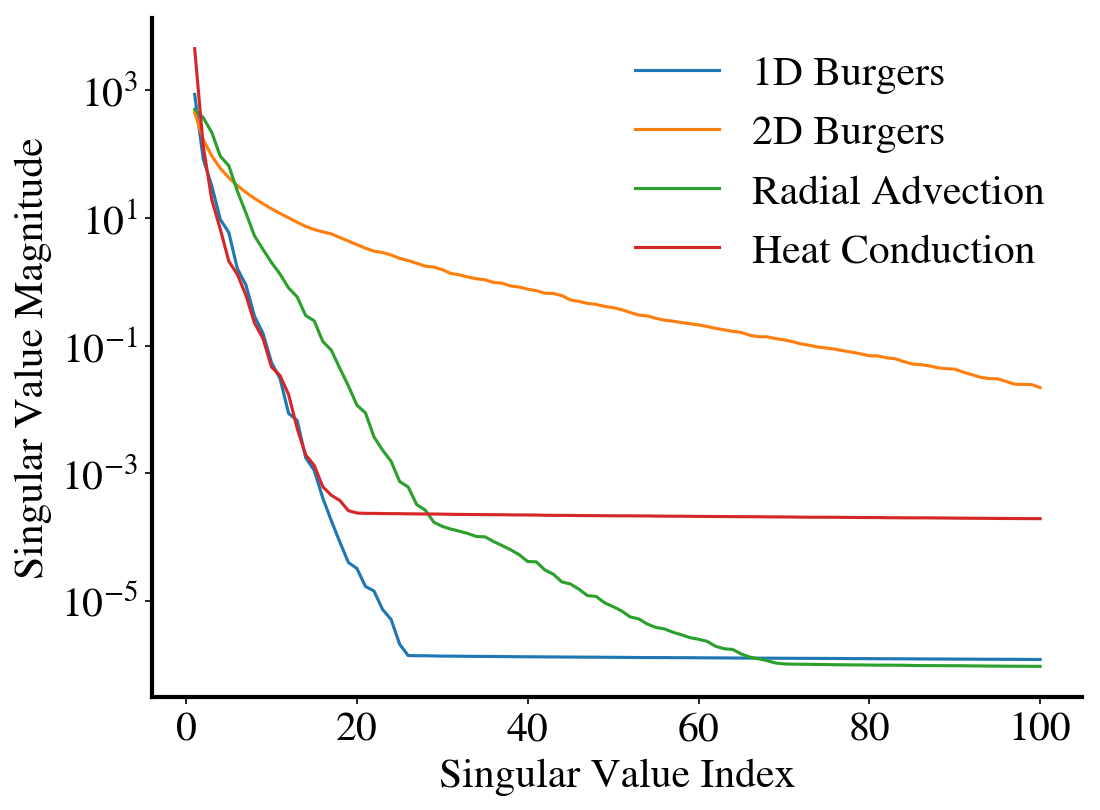}
    \caption{The decay of the singular values for all the problems considered in the paper. It shows that 2D Burgers and radial advection problems have slow singular value decay, indicating that nonlinear compression would work better than linear compression. On the other hand, 1D Burgers and heat conduction problems have fast singular value decay, indicating that linear compression is sufficient to achieve a good accuracy.}
    \label{fig: SV Decay}
\end{figure}

The second step of LaSDI is to compress the snapshot matrix $\snapshotMat$ either using linear or nonlinear compression techniques. The choice of either linear or nonlinear compression can be determined by the Kolmogorov $n$-width, which quantifies the optimal linear subspace. It is defined as:
\begin{equation}\label{eq:kolmogorov_nwidth}
    d_n(\mathcal{M}) := \inf_{\mathcal{L}_n} \sup_{f\in\mathcal{M}}\inf_{g\in\mathcal{L}_n} \| f-g \|,
\end{equation}
where $\mathcal{M}$ denotes the manifold of solutions over all time and parameters and $\mathcal{L}_n$ denotes all $n$-dimensional subspace. If a problem has a solution space whose Kolmogorov $n$-width decaying fast, then the the linear compression will provide an efficient subspace that can accurately approximates a true solution. However, a problem has a solution space with Kolmogorov $n$-width decaying slowly, then the linear compression will not be sufficient for good accuracy. Then, the nonlinear compression is necessary. The rate of decay in Kolmogorov $n$-width can be well indicated by the singular value decay (see Figure~\ref{fig: SV Decay}).
Section~\ref{sec:POD} describes a linear compression technique, i.e., proper orthogonal decomposition, which leads to \textbf{LaSDI-LS} where LS stands for ``Linear Subspace." Section~\ref{sec:auto-encoder} describes a nonlinear compression technique, i.e., auto-encoder, which leads to \textbf{LaSDI-NM} where NM stands for ``Nonlinear Manifold."

\subsubsection{Proper orthogonal decomposition: LaSDI-LS}\label{sec:POD}
We follow the method of snapshots first introduced by Sirovich \cite{sirovich1987turbulence}. The spatial basis from POD is an optimally compressed representation of $\snapshotMat$ in a sense that it minimizes the projection error, i.e., the difference between the original snapshot matrix and the projected one onto the subspace spanned by the basis, $\basis$:
  \begin{equation}\label{eq:POD}
    \begin{aligned}
      \basis &:= \underset{\basisDummy\in\RR{\nspacedof\times\nreducedspace},
      \basisDummy^T\basisDummy=\identityArg{\nreducedspace}}{\arg\min} & & \left \|\snapshotMat -
      \basisDummy\basisDummy^T{\snapshotMat} \right \|_F^2, 
    \end{aligned}
  \end{equation}
where $\|\cdot\|_F$ denotes the Frobenius norm. The solution of POD can be obtained by setting $\basis = \bmat{\leftSingularVecArg{1} & \cdots & \leftSingularVecArg{\nreducedspace}}$, $\nreducedspace < \nparam(\ntimedof+1)$, where $\leftSingularVecArg{\paramIndex}$ is $\paramIndex$th column vector of the left singular matrix, $\leftSingularMat$, of the following Singular Value Decomposition (SVD),
\begin{equation}\label{eq:SVD}
    \snapshotMat = \leftSingularMat \singularValueMat \rightSingularMat.
\end{equation}
Once the basis $\basis$ is built, the snapshot matrix $\snapshotMat\in\RR{\nspacedof\times(\ntimedof+1)\nparam}$ can be reduced to the generalized coordinate matrix, $\reducedSnapshotMat\in\RR{\nreducedspace\times(\ntimedof+1)\nparam}$, so called, the reduced snapshot matrix in the subspace spanned by column vectors of $\basis$, i.e., 
\begin{equation}\label{eq:reducedSnapshotMat}
    \reducedSnapshotMat := \basis^T\snapshotMat
\end{equation}
Naturally, we can define the reduced snapshot matrix, $\reducedSnapshotMatArg{\paramIndex}$, tailored for $\param_{\paramIndex}$ by extracting from $((\paramIndex-1)(\ntimedof+1)+1)$th to $\paramIndex(\ntimedof+1)$th column vectors of $\snapshotMat$, i.e., 
\begin{equation}\label{eq:reducedSnapshotMatArg}
    \reducedSnapshotMatArg{\paramIndex} := \bmat{\redsolapproxArgTwo{0}{\paramIndex} & \cdots & \redsolapproxArgTwo{\ntimedof}{\paramIndex}}, 
\end{equation}
where $\redsolapproxArgTwo{\timeIndex}{\paramIndex}\in\RR{\nreducedspace}$ is the reduced coordinates at $\timeIndex$th time step for $\param = \param_{\paramIndex}$. The matrix $\reducedSnapshotMatArg{\paramIndex}$ describes the latent space trajectory corresponding to $\param_{\paramIndex}$. For example, Step 2 of Figure~\ref{fig:procedureLaSDI} shows a graph of the latent space trajectory for the last parameter value, $\param_{\nparam}$. These reduced coordinates data will be used to train either global or local DI models to identify the dynamics of the latent space (see Section~\ref{sec:dynamicsIndentification}).

\subsubsection{Auto-encoder: LaSDI-NM} \label{sec:auto-encoder}
Auto-encoders act as a nonlinear analogue to POD. As illustrated in \cite{kim2021fast}, nonlinear subspace generated by auto-encoders outperform those generated by POD. In general, we train two neural networks, $\encoder: \RR{N_s} \rightarrow \RR{n_s}$ and $\decoder: \RR{n_s} \rightarrow \RR{N_s}$ to minimize 
\begin{equation}
\text{MSE}(\snapshotMat - \decoder(\encoder(\snapshotMat))
\end{equation}
where MSE denotes the mean-squared error. As above, we can define 
\begin{equation}
    \reducedSnapshotMat := \encoder(\snapshotMat)
\end{equation}
and can extract the snapshot $\reducedSnapshotMat_{k}$ by considering the $(k-1)(N_t+1) + 1)$th to $k(N_t+1)$th columns of $\reducedSnapshotMat$. 

The general architecture of $\encoder$ and $\decoder$ may vary. For the purposes of this paper, we will use a masked-shallow network as described in \cite{kim2021fast}. The use of this network architecture allows for universality across various PDEs and data shapes. The masking of the network allows for increased efficiency by not including full linear layers. However, the masking requires that we take care when constructing neural networks for higher dimensional simulations because the organization of spatial data must match the masking. Specifically, the architecture of both $\encoder$ and $\decoder$ consists of three layers, i.e., input, output and one hidden layer. The first layers are fully connected with nonlinear activation functions, whereas the final layer does not have nonlinear activation functions, i.e., fully linear. Then, the a masking operator is applied to make $\decoder$ sparse. We make the sparsity of the mask matrix similar to the one in the mass matrix to respect the sparsity induced by the underlying numerical scheme. Further details on the architecture of the neural networks can be found in \cite{kim2021fast}.

When training the auto-encoder, $n_s$ becomes a key hyper-parameter. Similar to POD, if $n_s$ is too small, the compression-decompression technique incurs a significant loss of data. Thus, $n_s$ cannot be chosen too small. In contrast to the POD, where increasing the basis size, improves results, this is not always the case with an auto-encoder. Rather, for larger $n_s$, significant improvement the error might not be seen. Further, when applying the DI techniques as described below, high-dimensional complex nonlinear systems are harder to approximate than those with simpler dynamics. Thus, we tune $n_s$ to be the smallest possible dimension while not compromising the accuracy of the reconstructed snapshot matrix. 

\subsection{Dynamics Identification} \label{sec:dynamicsIndentification}
To identify the latent space dynamics, we employ regression methods inspired by SINDy \cite{SINDy2016}. SINDy uses sparse regression tools, such as a sequential threshold least-squares, Lasso, Ridge, and Elastic Net. However, we do not require that our discovered dynamical system to be sparse in LaSDI. Instead, we want to generate the dynamical system that, when integrated, generates the most numerically accurate results. Because we do not require our solution to be sparse, we will refer this identification step to \textit{dynamics identification} (DI). 

In general, we aim to find a dynamical system $\dot{\redsolapprox}(t) = f(\redsolapprox(t))$ whose discrete solution best-matches the latent space discrete trajectory data, $\reducedSnapshotMat_{k}\in \RR{n_s\times (N_t+1)}$ in a least-squares sense. The goal is to identify a good $f(\redsolapprox(t))$. To formally describe the DI procedure, we take the transpose of the latent space trajectory data, i.e., ${\reducedSnapshotMat_{k}}^T$, to arrange the temporal evolution in a downward direction:
\begin{equation} \label{eq:temporaldownward}
   {\reducedSnapshotMat_{k}}^T = \bmat{\redsolapproxArg{0}^T \\ \vdots \\ \redsolapproxArg{N_t}^T} = \bmat{\redsolSymb_0(\timeSymb_0) & \hdots & \redsolSymb_{n_s}(\timeSymb_0) \\ \vdots & \ddots & \vdots \\ \redsolSymb_0(\timeSymb_{N_t}) & \hdots & \redsolSymb_{n_s}(\timeSymb_{N_t}) }
\end{equation}
Then we approximate the time derivative term with a finite difference and arrange the temporal evolution in a downward direction as well:
\begin{equation} \label{eq:dt_downward}
   {\dot{\reducedSnapshotMat}_{k}}^T = \bmat{\reddotsolapproxArg{0}^T \\ \vdots \\ \reddotsolapproxArg{N_t}^T} = \bmat{\dotredsolSymb_0(\timeSymb_0) & \hdots & \dotredsolSymb_{n_s}(\timeSymb_0) \\ \vdots & \ddots & \vdots \\ \dotredsolSymb_0(\timeSymb_{N_t}) & \hdots & \dotredsolSymb_{n_s}(\timeSymb_{N_t}) }
\end{equation}
We prescribe a library of functions, i.e., $\funcLibrary \left ({\reducedSnapshotMat_{k}}^T \right )\in\RR{(\ntimedof+1)\times\nlibrary}$, to approximate the right-hand side of the dynamical system, i.e., $f\left (\redsolapprox(t) \right )$. For example, $\funcLibrary \left ({\reducedSnapshotMat_{k}}^T \right )$ may take a form of constant, polynomial, trigonometric, and exponential functions:
\begin{equation} \label{eq:library}
  \funcLibrary \left ({\reducedSnapshotMat_{k}}^T \right ) = \bmat{1 & {\reducedSnapshotMat_{k}}^T & {\reducedSnapshotMat_{k,P_2}}^T & \hdots & {\reducedSnapshotMat_{k,P_{\ell}}}^T & \hdots & \sin({\reducedSnapshotMat_{k}}^T) & \cos({\reducedSnapshotMat_{k}}^T) & \hdots & \exp({\reducedSnapshotMat_{k}}^T) },
\end{equation}
where $\nlibrary$ denotes the number of columns, which is determined by the choice of a library of functions, and $P_{\ell}$ represents all the polynomials with order $\ell$. For the numerical examples in this paper, $\ell\leq 5$ are sufficient for accurate approximation of the dynamical system. For example, 
${\reducedSnapshotMat_{k,P_2}}^T$ and $\exp({\reducedSnapshotMat_{k}}^T)$ are defined as
\begin{equation} \label{eq:P2}
    {\reducedSnapshotMat_{k,P_2}}^T = \bmat{\redsolSymb_0^2(\timeSymb_0) & \redsolSymb_0(\timeSymb_0)\redsolSymb_1(\timeSymb_0) & \hdots & \redsolSymb_1^2(\timeSymb_0) & \hdots & \redsolSymb_{n_s}^2(\timeSymb_0) \\
    \vdots & \vdots & \ddots & \vdots & \ddots & \vdots \\
    \redsolSymb_0^2(\timeSymb_{N_t}) & \redsolSymb_0(\timeSymb_{N_t})\redsolSymb_1(\timeSymb_{N_t}) & \hdots & \redsolSymb_1^2(\timeSymb_{N_t}) & \hdots & \redsolSymb_{n_s}^2(\timeSymb_{N_t})
    }
\end{equation}
and 
\begin{equation} \label{eq:exo}
\exp({\reducedSnapshotMat_{k}}^T) = 
 \bmat{\exp(\redsolSymb_0(\timeSymb_0)) & \hdots &  \exp(\redsolSymb_{n_s}(\timeSymb_0)) \\
    \vdots & \ddots & \vdots  \\
    \exp(\redsolSymb_0(\timeSymb_{N_t})) & \hdots &  \exp(\redsolSymb_{n_s}(\timeSymb_{N_t}))}.
\end{equation}
The introduction of ${\reducedSnapshotMat_{k}}^T$, ${\dot{\reducedSnapshotMat}_{k}}^T$, and $\funcLibrary \left ({\reducedSnapshotMat_{k}}^T \right )$ in Eqs.~\eqref{eq:temporaldownward}, \eqref{eq:dt_downward}, and \eqref{eq:library} allows us to write the following regression problem:
\begin{equation}\label{eq:DI_single}
    \underset{\DIcoefficients_k\in\RR{\nlibrary\times\nreducedspace}}{\text{minimize}} \quad \left \| \dot{\reducedSnapshotMat}_{k}^T - \funcLibrary(\reducedSnapshotMat_k^T) \DIcoefficients_k \right \|^2_F
\end{equation}

So far, we have shown how to build local DI for a specific sample parameter, i.e., $\param_k$. However, this is not the only option. As a matter of fact, one can build a DI model for multiple parameter values. This observation leads to two types of DI: one, using all the training data to inform the underlying dynamical system; two, using only training data close to a query parameter point to model the dynamics. The former we refer to as \textit{Global DI} and the latter, \textit{Local DI}. Note that both global and Local DIs can be considered as \textbf{region-based} because they set a unified DI for each region. 
Alternatively, we can consider \textbf{point-wise} DIs, where the coefficients, $\DIcoefficients_k$ of each DI model built for a single training point, i.e., $\param_k$, by Eq.~\eqref{eq:DI_single} can be interpolated for a new query point, $\param^*$. We refer the detailed explanation of such interpolation to \textit{interpolated DI}. 

\subsubsection{Region-based global DI}\label{sec:globalDI}
Global DI requires a straight-forward implementation. We include all the training snapshots in the identification of the latent space dynamics. For this, the discovered dynamics come from modifying \eqref{eq:DI_single} to 
\begin{equation}\label{eq:globalDI}
    \underset{\DIcoefficients\in\RR{\nlibrary\times\nreducedspace}}{\text{minimize}} \quad \sum_{k=1}^{\nparam}\left \| \dot{\reducedSnapshotMat}_k^T - \funcLibrary(\reducedSnapshotMat_k^T) \DIcoefficients \right \|^2_F
\end{equation}
This gives one universal DI model that covers the whole parameter space, which might be a good model if the latent space dynamics do not drastically change based on parameter values. This might be valid for a small parameter space. However, for a larger parameter space, the latent space dynamics might change drastically. Therefore, we introduce local DI.

\subsubsection{Region-based local DI}\label{sec:localDI}
\begin{figure}
\centering
    \includegraphics[width = .7\textwidth]{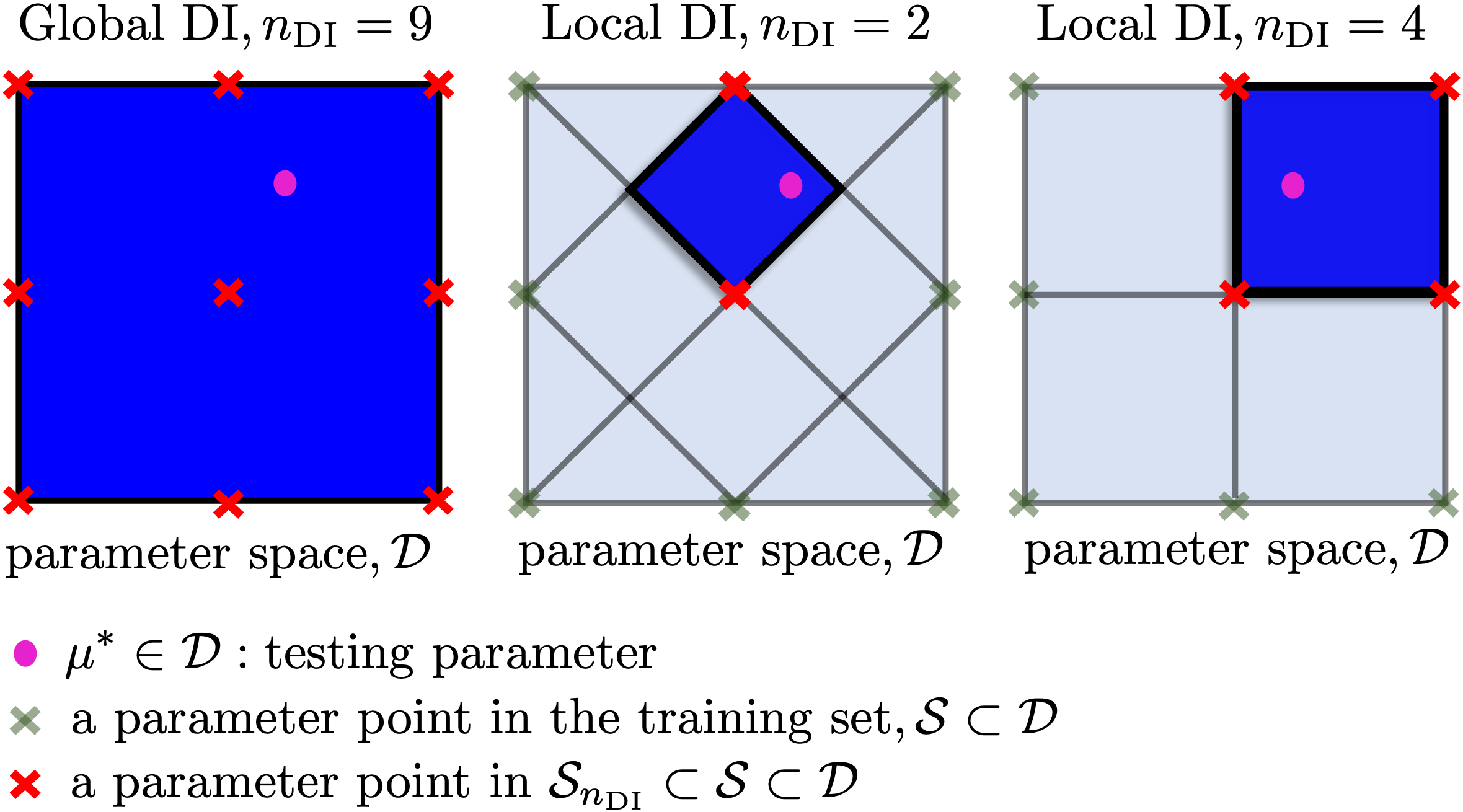}
    \caption{An example of a 2-dimensional parameter space with nine training values uniformly spaced throughout the space. For an arbitrary parameter point, $\param^*\in\paramDomain$, we illustrate three different latent space dynamics identification regions, $\region{\localDIindex}$, colored by blue, for $\numDI=9, 2$, and $4$. The case of $\numDI=9$ is the global DI, while the cases of $\numDI=2$ and $4$ are local DIs.}
    \label{fig: method local global}
\end{figure}

For a large parameter space, we should not expect that a single dynamical system would govern all the latent space trajectories from the whole training data. However, we might expect that the dynamics in the latent space locally behave. This observation allows us to build a local DI which provides a good model for a local region. The local DI takes training simulation data from $\numDI$ neighboring points. Therefore, the least-squares regression problem for $\localDIindex$th local DI becomes
\begin{equation}\label{eq:localDI}
    \underset{\DIcoefficients_{\localDIindex}\in\RR{\nlibrary\times\nreducedspace}}{\text{minimize}} \quad \sum_{k\in\indexSet(\localDIindex)}\left \| \dot{\reducedSnapshotMat}_k^T - \funcLibrary(\reducedSnapshotMat_k^T) \DIcoefficients_{\localDIindex} \right \|^2_F,
\end{equation}
where $\setLocalDI(\localDIindex)\in\trainingSet$ denotes a set of $\numDI$ training parameters that determine $\localDIindex$th local parameter region and its corresponding set of indices is denoted as $\indexSet(\localDIindex)$. The $\localDIindex$th local parameter region, $\region{\localDIindex}$, is defined by the set of points that are closest to parameters in $\setLocalDI(\localDIindex)$ with some distance measure, i.e.,
\begin{equation}
    \region{\localDIindex} \equiv \{ \param^*\in\paramDomain \mid \dist{\param^*}{\param_j} \leq \dist{\param^*}{\param_k}, \param_j\in\setLocalDI(\localDIindex), \param_k\in\trainingSet\setminus\setLocalDI(\localDIindex) \},
\end{equation}
where $\dist{\cdot}{\cdot}$ denotes a distance between two points in $\paramDomain$. We use Euclidean distance that is defined as
\begin{equation}\label{eq:nn}
    d(\param_i, \param_j) := ||\param_i-\param_j||_2,
\end{equation}
but other distance functions can be used as well. Note that $\numDI$ is a parameter that can be tuned. If $\numDI=\nparam$, we recover the global DI. Therefore, the global DI can be considered as a special case of local DIs. We illustrate the various examples of global and local DI in a 2D parameter space and their local parameter region which are relevant to a specific $\param^*$ denoted as pink dots in Figure \ref{fig: method local global}.

\subsubsection{Point-wise interpolated DI}\label{sec:interpolatedDI}
 For $\param^*\in\region{\localDIindex}$, the interpolated DI interpolates coefficients, $\DIcoefficients_k$, for $k\in\indexSet(\localDIindex)$ if $\param^*\in\region{\localDIindex}$, where
 $\DIcoefficients_k$ are obtained by solving Eq.~\eqref{eq:DI_single}, to obtain the interpolated coefficients, $\DIcoefficients^*$, which will be used for the LaSDI prediction in Eq.~\eqref{eq:LaSDI_ODE}. Many interpolation techniques can be utilized here, but we use two techniques: radial basis functions (RBFs) and linear bi-variate splines. In each case, we interpolate element-wise in $\DIcoefficients^*$. For the $\nlibrary\cdot \nreducedspace$ elements in $\DIcoefficients^*$, the corresponding interpolation functions are defined below.
 
The RBF interpolation function for each element of $\DIcoefficients^*$ for $\param^*\in\region{\localDIindex}$, is defined as
  \begin{equation}\label{eq:RBF}
     \rbfSymb(\param^*; \{\param_k\}_{k\in \indexSet(\localDIindex)}) = \sum_{k=1}^{\indexSet(\localDIindex)} \rbfcoefk \ \rbfSymb_k \left (\dist{\param^*}{\param_k} \right )
 \end{equation}

 where $\rbfSymb_k$ are radial basis functions; $\rbfcoefk\in\RR{}$ denotes $k$th RBF interpolation coefficient, which is obtained by solving a least-squares problem; $\dist{\cdot}{\cdot}$ is the Euclidean distance as defined in Eq.~\eqref{eq:nn}; and $\param_k\in\region{r} \cap \trainingSet$. Many radial basis functions are available. For our numerical experiments, we use the multiquadric radial basis functions defined as 
 \begin{equation}
     \rbfSymb_k (\param;\{\param_k\}_{k\in \indexSet(\localDIindex)}) = \sqrt{\frac{\dist{\param^*}{\param_k}^2}{\epsilon^2} + 1},
 \end{equation}
 where $\epsilon$ approximates average distance between $\{\param_k\}_{k\in\indexSet(\localDIindex)}$.

 The linear bi-variate splines approximates $\DIcoefficients^*$ using a linear combination of the $\DIcoefficients_k$ based on $\dist{\param^*}{\param_k}$ with no smoothing. If $\param_k$ falls on a uniform grid, as in our numerical experiments, the bi-linear interpolation function for each entry of $\DIcoefficients^*$ can be defined in rectangular regions. For example, in two dimensional parameter space, this is defined by the four corner points in a specific region (e.g., $\region{r}$), $\param_0 = \bmat{\paramSymb_{0,0} & \paramSymb_{0,1}}^T$, $\param_1 = \bmat{\paramSymb_{1,0} & \paramSymb_{0,1}}^T$,  $\param_2 = \bmat{\paramSymb_{0,0} & \paramSymb_{1,1}}^T$, $\param_3 = \bmat{\paramSymb_{1,0} & \paramSymb_{1,1}}^T$ and the corresponding $\DIcoefficients_k$s. For $\param^* = \bmat{\paramSymb^*_0 & \paramSymb^*_1}^T\in \region{r}\setminus\partial\region{r}$, each entry of $\DIcoefficients^*$ can be obtained by
\begin{equation}\label{eq:bilinear}
    \rbfSymb(\param^*;\{\param_k\}_{k\in \indexSet(\localDIindex)}) = \frac{1}{(\paramSymb_{1,0}-\paramSymb_{0,0})(\paramSymb_{1,1}-\paramSymb_{0,1})}[\paramSymb_{1,0}-\paramSymb^*_0\;\; \paramSymb^*_0-\paramSymb_{0,0}]\begin{bmatrix} \DIcoefficients_{0} & \DIcoefficients_{1} \\ \DIcoefficients_{2} & \DIcoefficients_{3}\end{bmatrix}\begin{bmatrix} \paramSymb_{1,1}-\paramSymb^*_1 \\ \paramSymb^*_1-\paramSymb_{0,1}\end{bmatrix}
 \end{equation}

\subsubsection{Notes on DIs}\label{sec:notesOnDI}
 In general, finding an accurate and numerically stable approximation can be difficult for unknown nonlinear problems. Because of this, the terms included in the regression becomes a hyper-parameter that must be tuned. As with tuning $n_s$, we strive to use the least number of terms for best accuracy. Therefore, we now discuss the procedure used to obtain the best dynamics approximation given a set of trajectories in a latent space:

\begin{enumerate}
    \item \textbf{Start simple}: start with polynomial order one. This generates the linear dynamical system that best approximates the latent space dynamics.
    \item \textbf{Verification}: Visually verify that the new trajectories approximate the training data.
    \item If a better fit is required, we include the following modifications in order of importance:
    \begin{itemize}
        \item \textbf{Rescale}: rescale the data so that $\max\left(\left|\reducedSnapshotMat_{k}\right|\right)\leq 1$. 
        \item \textbf{Enrich}: enrich the library with more terms, such as polynomials with a higher order, mixed polynomial terms (e.g., $\redsolSymb_0\redsolSymb_1$), trigonometric and exponential terms. 
    \end{itemize}
\end{enumerate}

Solving one of the least-squares problems, i.e., Eqs.~\eqref{eq:DI_single}, \eqref{eq:globalDI}, and \eqref{eq:localDI}, or utilizing point-wise approach in Section~\ref{sec:interpolatedDI} gives a latent space dynamical system characterized by the form of a system of ordinary differential equations (ODEs):
\begin{equation} \label{eq:LaSDI_ODE}
   \dot{\redsolapprox}(t) = \rhsArg{\numDI,\localDIindex}(\redsolapprox(t);\DIcoefficients),
\end{equation}
where right-hand-side, $\rhsArg{\numDI,\localDIindex}:\RR{n_s}\times\RR{\nlibrary}\rightarrow\RR{n_s}$, is determined by a library of functions specified in Eq.~\eqref{eq:library} and the coefficients identified from the least-squares problems or interpolations.

\subsection{Prediction}\label{sec:prediction}
After local or global or interpolated DIs are set in the dynamics identification step, one can solve them to predict new latent space dynamics trajectories corresponding to a new parameter $\param^*\in\paramDomain$. If the global DI is used, there is no ambiguity. However, if a local DIs are used, then one needs to figure out which local DI needs to be solved. It can be determined by checking  $\param^*\in\region{\localDIindex}$. However, pre-defining $\region{\localDIindex}$ over the whole parameter space, $\paramDomain$, can be a daunting process for non-uniform training points. Therefore, the interpolated DIs might be more practical for non-uniform training points.

\subsubsection{Latent space dynamics solve and reconstruction}\label{sec:latentspace_solve}
Once an appropriate DI is recognized for a given $\param^*$, an appropriate initial condition $\redsolapproxArg{0}^*$ is required. The latent space initial condition corresponding to $\param^*$ is obtained by applying compression step to the full order model initial condition $\solArg{0}^*$ corresponding to $\param^*$ as described in Section~\ref{sec:compression}. For the linear compression, 
\begin{equation} \label{eq:linear_init}
    \redsolapproxArg{0}^* \equiv \basis^T\solArg{0}^*
\end{equation}
and for the nonlinear compression,
\begin{equation} \label{eq:nonlinear_init}
    \redsolapproxArg{0}^* \equiv \encoder(\solArg{0}^*)
\end{equation}
are used. Using this initial condition, Eq.~\eqref{eq:LaSDI_ODE} is solved to predict the latent space dynamics trajectory, $\redsolapproxArg{n}^*$ using the Dormand-Prince pair of formulas \cite{dormand1980family}. At each iteration, the local error is controlled using a fourth-order method. Using local extrapolation, the step is then taken using a fifth-order formulation. Then, the approximated full state trajectories, $\solapproxArg{n}^*$ are restored by $\solapproxArg{n}^* \equiv \basis\redsolapproxArg{n}^*$ if the linear compression were used and by $\solapproxArg{n}^* \equiv \decoder(\redsolapproxArg{n}^*)$ if the autoencoder is used as a nonlinear compression.

The auto-encoder training and dynamics identification steps can be simultaneously done as in \cite{champion2019data, he2022glasdi} to obtain simpler latent space dynamics. However, this paper demonstrates that the auto-encoder training and dynamics identification steps can be separated.

\section{Numerical Results}\label{sec:results}

\label{sec: Results}

We demonstrate the performance of LaSDI for four different numerical examples, i.e., 1D and 2D Burgers equations, heat conduction, and radial advection problems. The governing equations, initial conditions, boundary conditions, parameter space, and domain are specified in Table \ref{tab:examples}. 

\begin{table}[ht]
    \centering
    \def\arraystretch{1.4}
    \begin{tabular}{|c|c|c|}
        \hline
        \textbf{Equation} & \textbf{Initial Condition} & \textbf{Domain/Boundary} \\ \hline
        \begin{tabular}{c}
        \underline{1D Burgers} \\
        $u_t = -u  u_x$ \end{tabular} & $\begin{array}{c} u(0,x;a,w) = a\cdot \exp(-x^2/w) \\ a\in [0.7, 0.9] \\ w\in [0.9,1.1] \end{array} $ & $\begin{array}{c} \Omega =[-3,3]\\ t\in [0,1]\\ u(-3,t) = u(3,t) = 0 \end{array}$\\ 
        \hline
        \begin{tabular}{c}
        \underline{2D Burgers}\\
        $\mathbf{u}_t = - \mathbf{u}\cdot \nabla \mathbf{u} + \frac{1}{10000}\Delta \mathbf{u}$ \end{tabular} & $\begin{array}{c} \mathbf{u}(0,\mathbf{x};a,w) = \begin{bmatrix} a\cdot \exp(-||\mathbf{x}||_2^2/w) \\ a\cdot \exp(-||\mathbf{x}||_2^2/w) \end{bmatrix} \\ a\in [0.7, 0.9]\\w\in [0.9,1.1]\end{array}$ & $\begin{array}{c}\Omega =[-3,3]\times [-3,3] \\ t\in [0,2] \\  \mathbf{u}(\mathbf{x},t) = 0\text{ on } \partial\Omega \end{array}$\\
        \hline
        \begin{tabular}{c}
        \underline{Heat Conduction}\\
        $u_t= \nabla\cdot (1+u)\nabla u$ \end{tabular} & $\begin{array}{c} u(0,
        \mathbf{x}; \omega, a) = a\sin(\omega(x_1+x_2))+a \\ \omega \in [0.2,5.0]\\a\in [1.8,2.2]\end{array}$ & $\begin{array}{c}\Omega = [0,1]\times [0,1] \\ t\in [0,1] \\  \frac{\partial u}{\partial n} = 0\text{ on } \partial \Omega  \end{array}$\\
        \hline
        \begin{tabular}{c}
        \underline{Radial Advection}\\
        $u_t =- \mathbf{v}\cdot \nabla u$ \\ $\mathbf{v} = \frac{\pi}{2}d\begin{bmatrix} x_2 \\ -x_1 \end{bmatrix}$ \\ $
        d =((1-x_0^2)(1-x_1^2))^2$
        \end{tabular} & $\begin{array}{c} u(0,\mathbf{x};\omega) = \sin(\pi \omega x_1)\sin(\pi\omega x_2) \\ \omega \in [0.6,1.0] \text{ or } \omega \in [0.6,1.4]\end{array}$ & $\begin{array}{c}\Omega = [-1,1]\times [-1,1] \\ t\in [0,3] \\  u(\mathbf{x},t) = 0\text{ on } \partial \Omega \\  \end{array}$\\
        \hline
    \end{tabular}
    \caption{The four PDEs used to generate FOMs for the application of LaSDI.}
    \label{tab:examples}
\end{table}

For both 1D and 2D Burgers problems, a uniform space discretization (i.e., $dx = 6/1000$ for 1D and $dx = dy = 1/10$ for 2D Burgers problems) is used for the full order model solve. The first order spatial derivative term is approximated by the backward difference scheme, while the diffusion terms are approximated by the central difference scheme. This generates a semi-discretized PDE characterized by a system of nonlinear ordinary differential equations specified in Eq.~\eqref{eq:dynamics}. We integrate each of these systems over a uniform time-step (i.e., $\Delta t = 1/1000$ for 1D and $\Delta t = 2/1500$ for 2D Burgers problems) using the implicit backward Euler time integrator, i.e., $\solArg{n} - \solArg{n-1} = \Delta t \rhsArg{n}$.

For the full order model of the nonlinear heat conduction problem, we use the second order finite element to discretize the spatial domain. The spatial discretization starts with $64\times64$ uniform squares, which are divided into triangular elements, resulting in 8192 triangular elements in total. A uniform time step of $\Delta t = .01$ is used with the SDIRK3 implicit L-stable time integrator (i.e., see Eq.~(229) of \cite{kennedy2016diagonally}). For the temperature field dependent conductivity coefficient, we linearize the problem by using the temperature field $\sol$ from the previous time step.\footnote{The source code of the full order model for the nonlinear heat conduction problem can be found at https://github.com/mfem/mfem/blob/master/examples/ex16.cpp}

For the full order model of the radial advection problem,  we use the third order discontinuous finite element to discretize the spatial domain. The spatial discretization is dictated by a square-mesh with periodic boundary conditions. We use $24\times24$ grid of square finite elements. The finite element data is interpolated to generate a $64\times64$ uniform grid across the spatial domain. We use a uniform time step of $\Delta t = .0025$ with the RK4 explicit time integrator.\footnote{The source code of the full order model for the radial advection problem can be found at https://github.com/mfem/mfem/blob/master/examples/ex9.cpp}

These full order models are used to generate training data. Table~\ref{tab:training} shows testing parameter set, $\testingSet$, as well as all the training samples, $\trainingSet$, for four problems. Note that the testing points include the training points. The accuracy is measured by the maximum relative error, $\relerr:\paramDomain\rightarrow\RRplus{}$, at each testing parameter point, which is defined as 
\begin{equation}
\label{eq: error}
    \relerr(\param^*) = \max_{n\innat{\ntimedof}} \left(\frac{||\solapproxArg{n}(\param^*)-\solArg{n}(\param^*)||_2}{||\solArg{n}(\param^*)||_2}\right),
\end{equation}
where $\param^*$ is a testing parameter point. 

\begin{table}[ht]
    \centering
    \def\arraystretch{1.5}
    \begin{tabular}{|c|c|c|}
        \hline
        \textbf{Example} & \textbf{Training set}, $\trainingSet$ & \textbf{Testing set}, $\testingSet$ \\ \hline
        1D Burgers & $\begin{array}{c} a\in \{0.7, 0.9\};\; w\in \{0.9,1.1\}, \\ a\in \{0.7, 0.8, 0.9\};\; w\in \{0.9,1.0,1.1\}, \\ \text{ or } \\ a\in \{0.7,0.75, \dots, 0.9\}; \\ w\in \{0.9,0.95, \dots, 1.1\}\end{array}$ & $ \begin{array}{c} a\in \{0.7,0.71,\dots, 0.9\} \\ w \in \{0.9, 0.91, \dots 1.1\} \end{array}$\\ 
        \hline
        2D Burgers & $\begin{array}{c}  a\in \{0.7,0.75, \dots, 0.9\}; \\ w\in \{0.9,0.95, \dots, 1.1\}\end{array}$ & $ \begin{array}{c} a\in \{0.7,0.72,\dots, 0.9\} \\ w \in \{0.9, 0.92, \dots 1.1\} \end{array}$\\
        \hline
        Heat Conduction & $\begin{array}{c} \omega \in \{0.2,0.8,\dots, 5.0\}; \\ a\in \{1.8,2.0,2.2\}\end{array}$ & $ \begin{array}{c} a\in \{0.2,0.24,\dots, 5.0\} \\ w \in \{1.8, 1.81, \dots 2.2\} \end{array}$\\
        \hline
        Radial Advection & $\begin{array}{c}\omega \in \{0.6,0.7,\dots, 1.0\}, \\ \omega \in \{0.6,0.65,\dots, 1.0\}, \\ \omega \in \{0.6,0.7,\dots, 1.4\}, \\ \text{ or } \\ \omega \in \{0.6,0.65,\dots, 1.4\}
        \end{array}$ & $ \begin{array}{c} \omega\in \{0.6,0.61,\dots, 0.1.0\} \\ \omega \in \{0.6, 0.61, \dots 1.4\} \end{array}$\\
        \hline
    \end{tabular}
    \caption{The specific parameter values used in training and testing LaSDI throughout the paper. Each subsequent figure makes clear the number of training values used for training the data-compression technique as well as the latent space DI.}
    \label{tab:training}
\end{table}

The computational cost is measured in terms of the wallclock time. Specifically, timing is obtained by performing calculations on a IBM Power9 @ 3.50 GHz and DDR4 Memory @ 1866 MT/s. The auto-encoders are trained on a NVIDIA V100 (Volta) GPU with 3168 NVIDIA CUDA Cores and 64 GB GDDR5 GPU Memory using PyTorch [55] which is the open source machine learning frame work.

\subsection{1D Burgers}\label{sec:1dB}

\begin{figure}
    \centering
        \includegraphics[width=.5\textwidth]{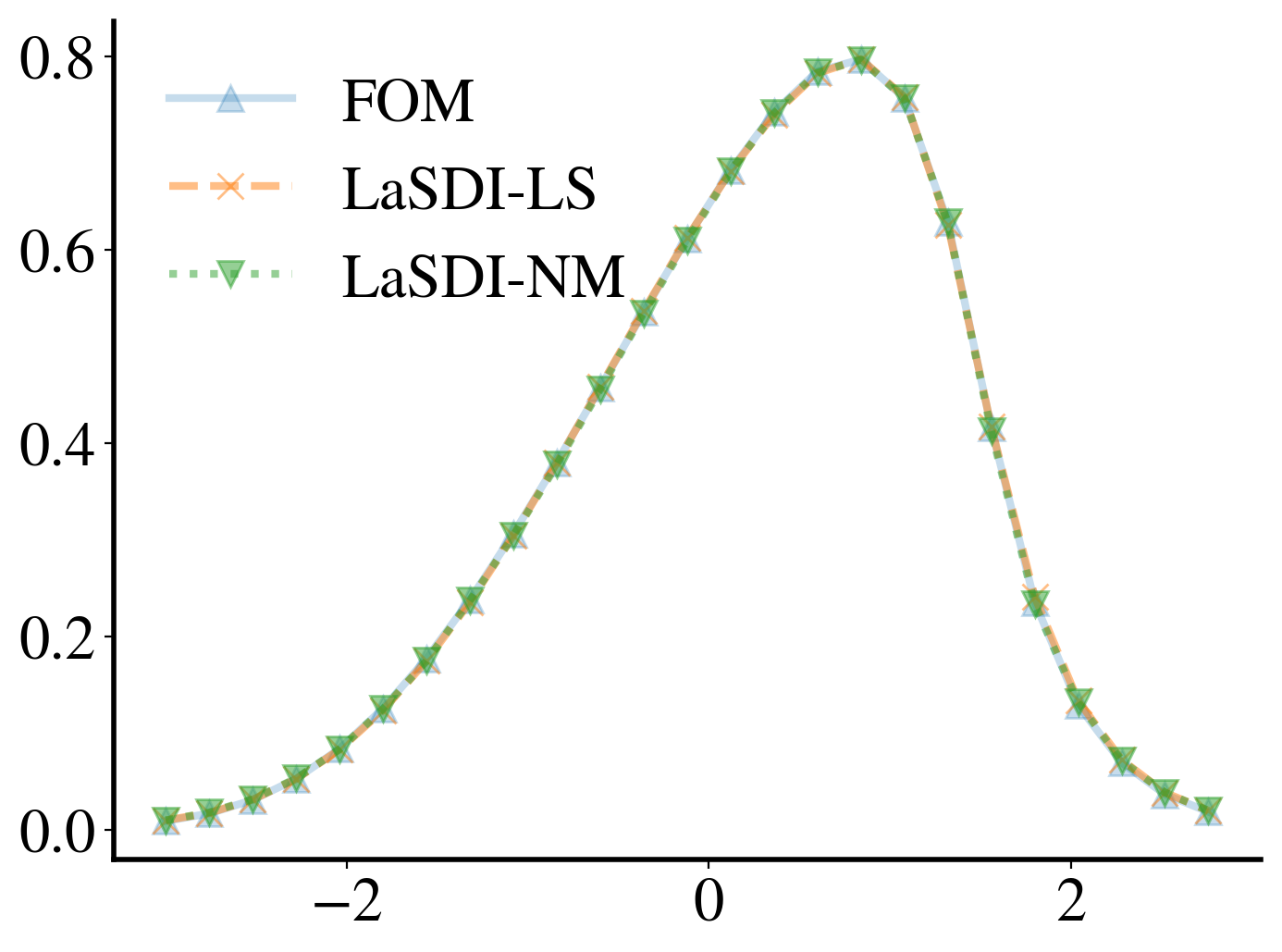}
        \caption{LaSDI ROMs at the final time step for the 1D Burgers problem generated using either the POD or nonlinear data-compression techniques for $\mu^* = (0.8,1.01)$. In both instances, we use linear dynamical systems in global DI. LaSDI-NM and LaSDI-LS use a latent-space dimension of four and five respectively.}
        \label{fig:1D single}
\end{figure}

\begin{figure}
    \centering
    \includegraphics[width=\linewidth]{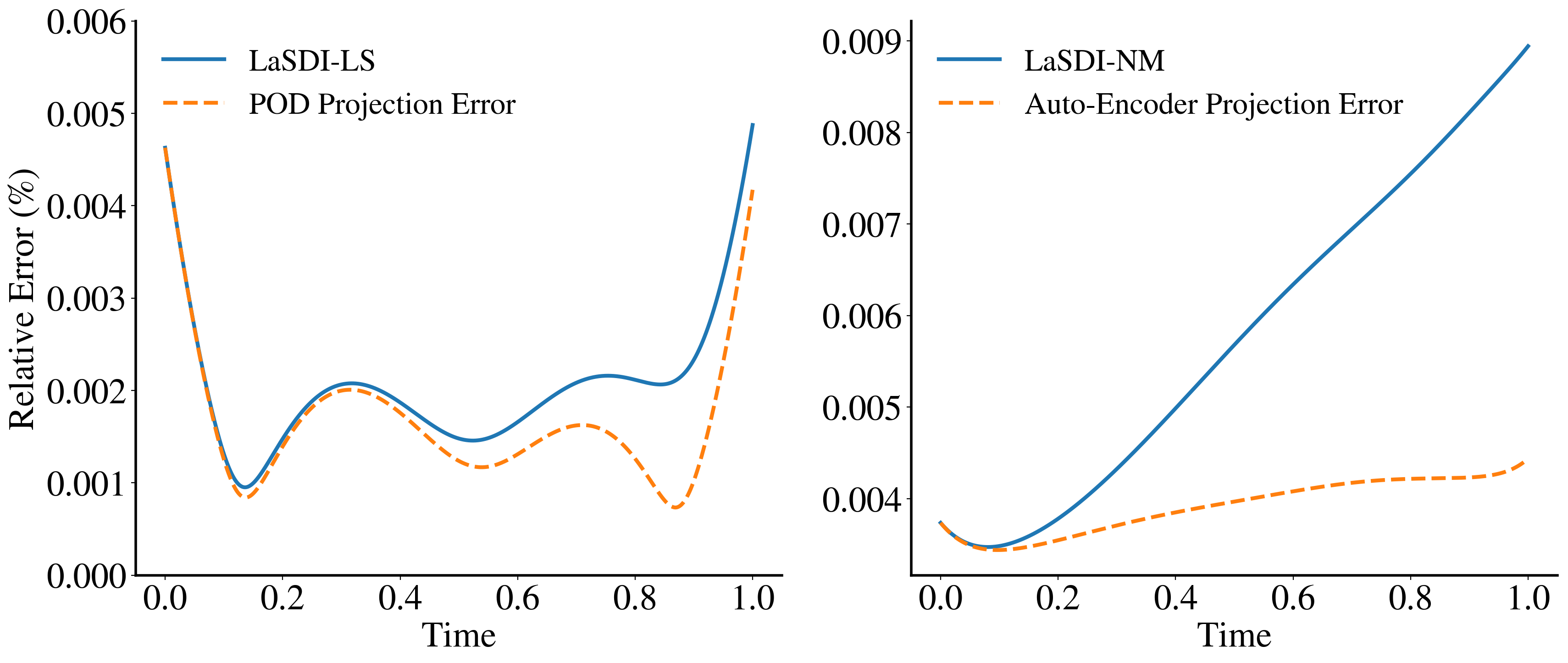}
\caption{The illustration of relative errors of LaSDI models being bounded below by the projection error of each compression technique, i.e., proper orthogonal decomposition and auto-encoder. The relative errors are computed from the simulation illustrated in Figure~\ref{fig:1D single}. The projection errors are computed by compressing and de-compressing the corresponding full order model solutions.}
\label{fig: 1db errors}
\end{figure}

First, for both LaSDI-LS and LaSDI-NM, we approximate the dynamics with linear dynamical systems in global DI. In Figure \ref{fig:1D single}, the last time step solutions for $\mu^* = (0.8, 1.01)$ are compared with the corresponding FOM solution, showing that they are almost identical. This shows that LaSDIs are able to accurately predict solutions. Figure \ref{fig: 1db errors} illustrates the relative error of a LaSDI model is bounded below by the projection error of the compression technique. 

\begin{figure}
    \centering
    \begin{subfigure}[b]{.45\linewidth}        \includegraphics[width=\textwidth]{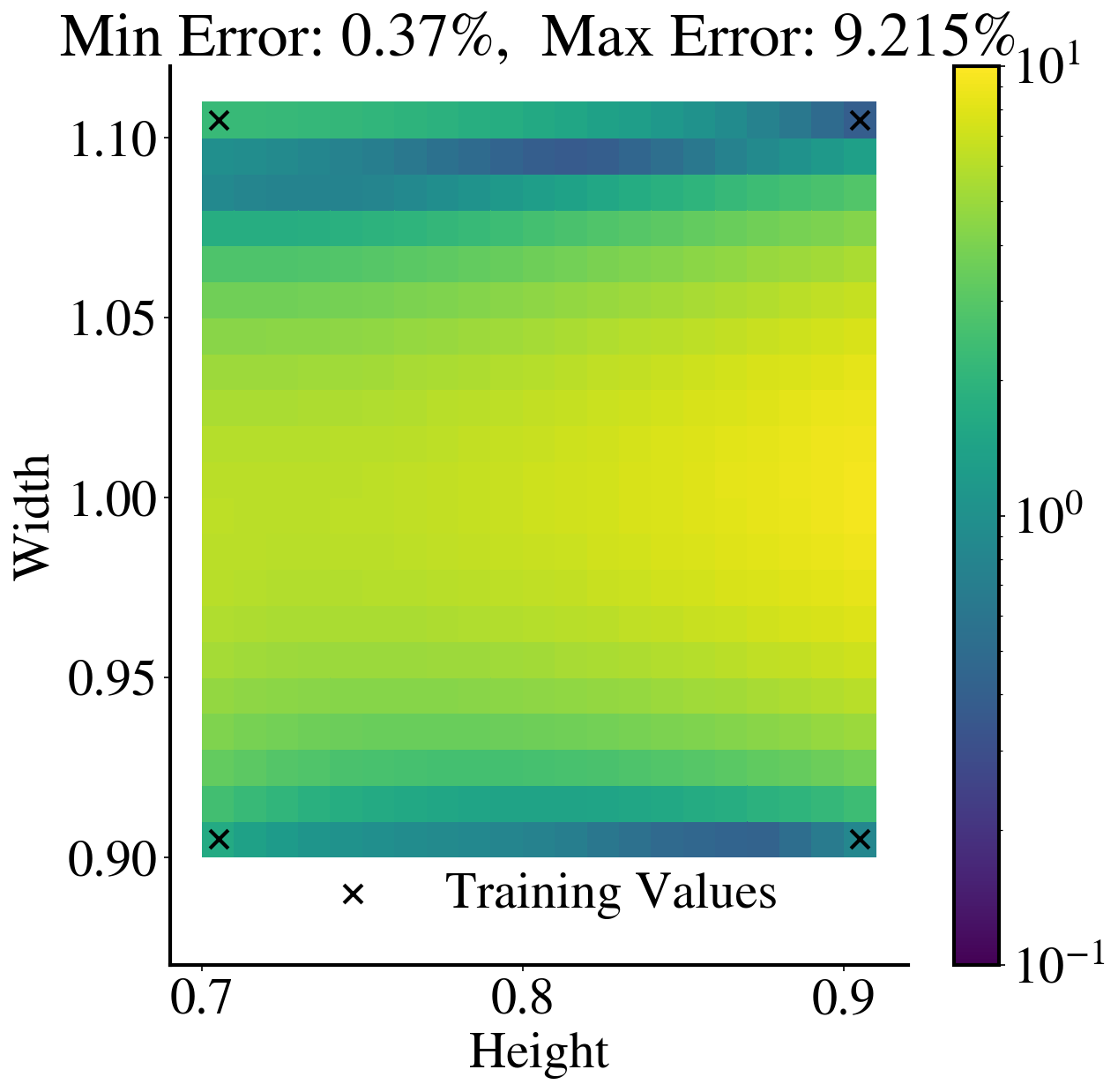}
        \caption{LaSDI-NM with latent-space dimension of 4.}
        \label{fig: 1db heat 4 NM}
    \end{subfigure}
    \begin{subfigure}[b]{.45\linewidth}
        \includegraphics[width=\textwidth]{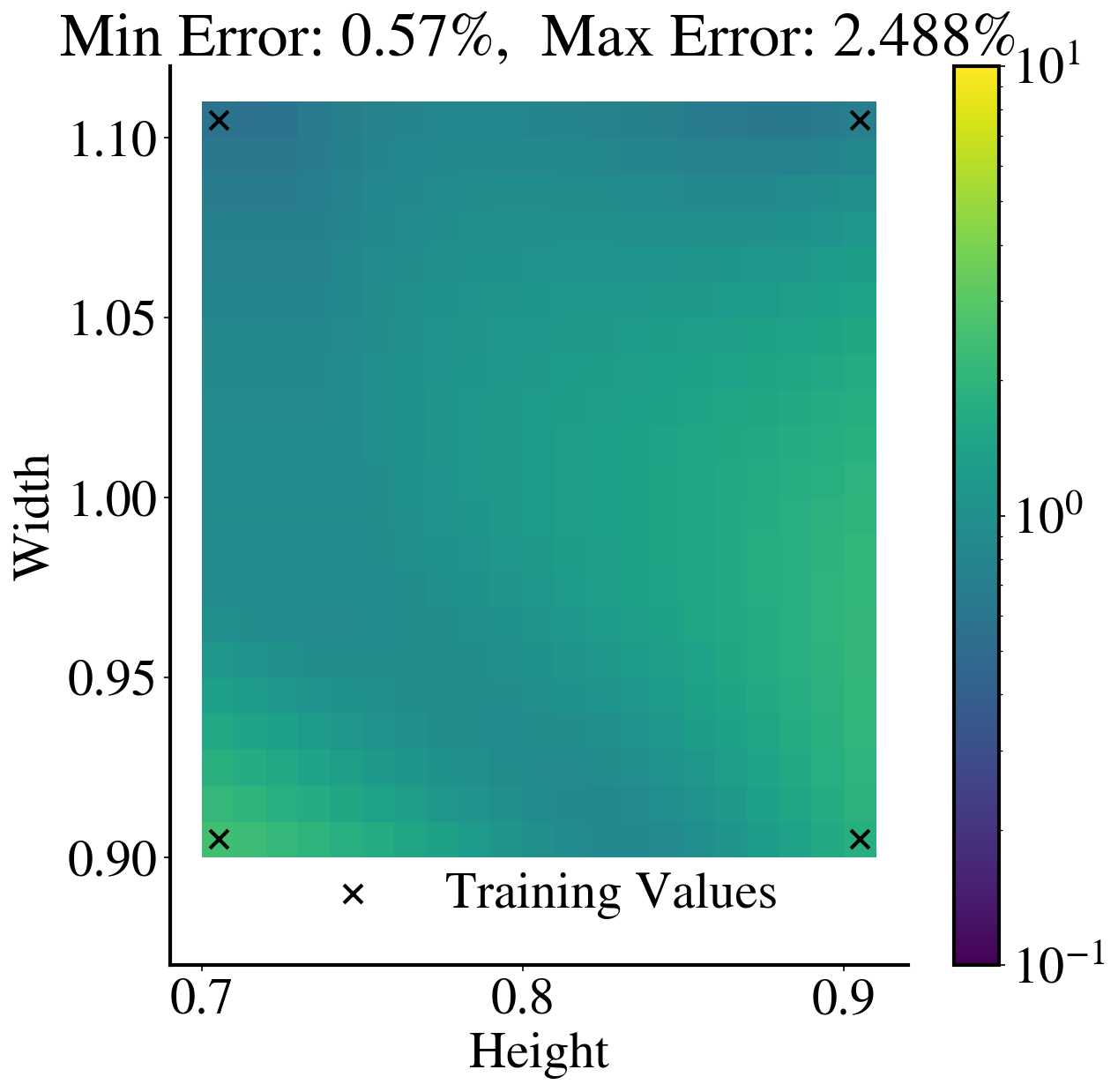}
        \caption{LaSDI-LS with latent-space dimension of 5.}
        \label{fig: 1db heat 4 ls}
    \end{subfigure}
\caption{Heat map of maximum relative errors on 1D Burgers parameter space both for LaSDI-NM and LaSDI-LS. Both LaSDI-NM and LaSDI-LS takes Global DI. Four training points are used to generate simulation data and global DI is used to predict the solution at each testing point. Both maximum and minimum relative errors are indicated.}
\label{fig: 1db heat 4}
\end{figure}

To see the accuracy over the whole parameter domain, Figure~\ref{fig: 1db heat 4} shows two heat maps of the maximum relative errors for each parameter case: one for LaSDI-NM and the other for LaSDI-LS. Both LaSDIs are trained using four training points. For this particular experiments, LaSDI-LS outperforms LaSDI-NM in terms of accuracy, i.e., The maximum relative error for LaSDI-NM is $9.2\%$, while the maximum relative error for LaSDI-LS is $2.5\%$.

\begin{figure}
    \centering
    \begin{subfigure}[b]{\linewidth}
        \centering
        \includegraphics[width = \linewidth]{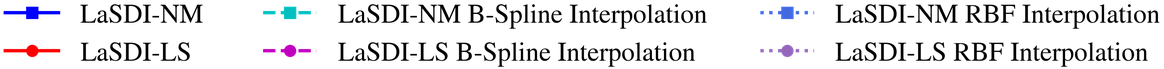}
    \end{subfigure}
    \begin{subfigure}[b]{.45\linewidth}
        \centering
        \includegraphics[width = \linewidth]{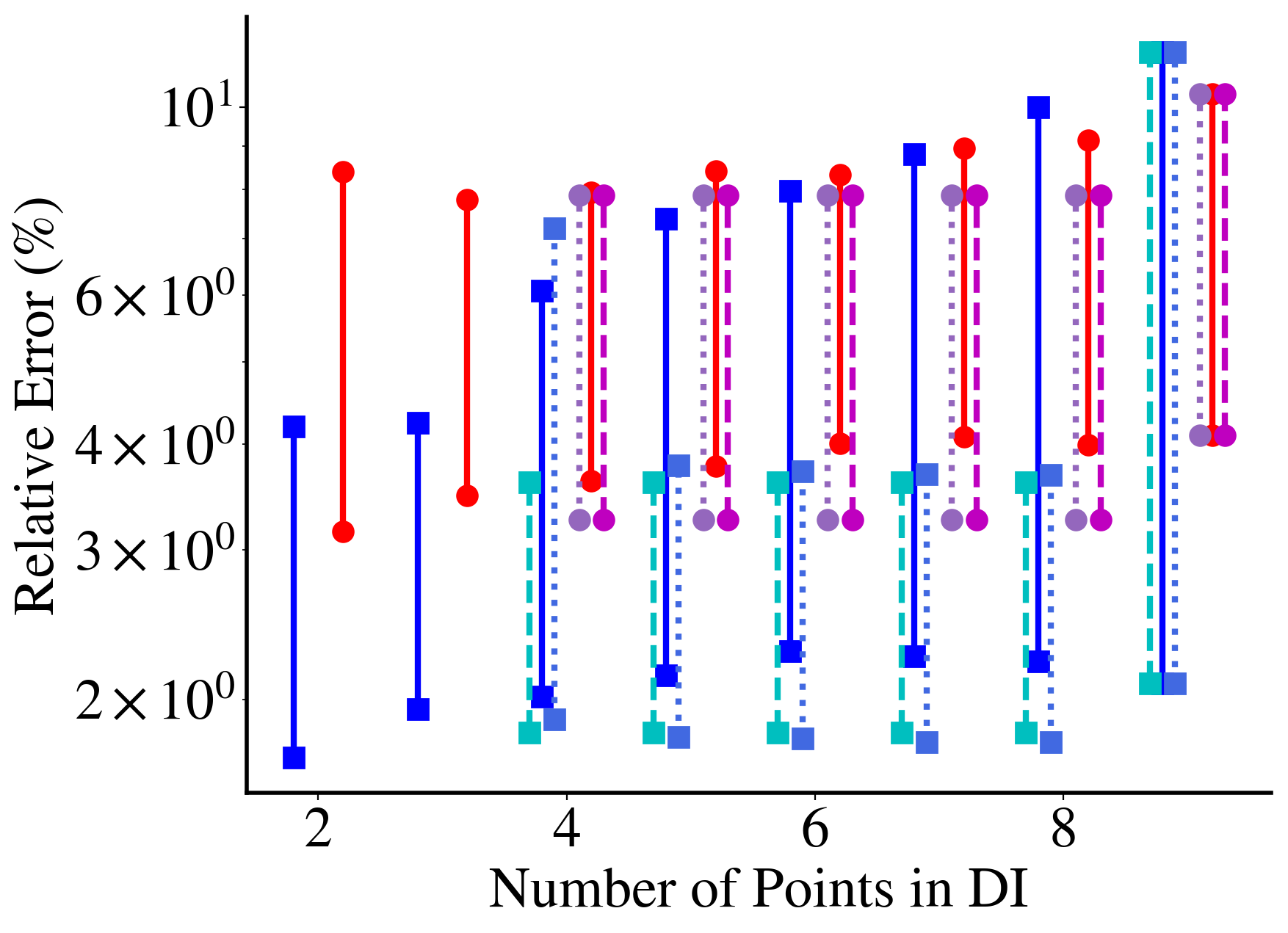}
        \caption{Error for degree 0 dynamical system in DI.}
    \end{subfigure}
        \begin{subfigure}[b]{.45\linewidth}
        \centering
        \includegraphics[width = \linewidth]{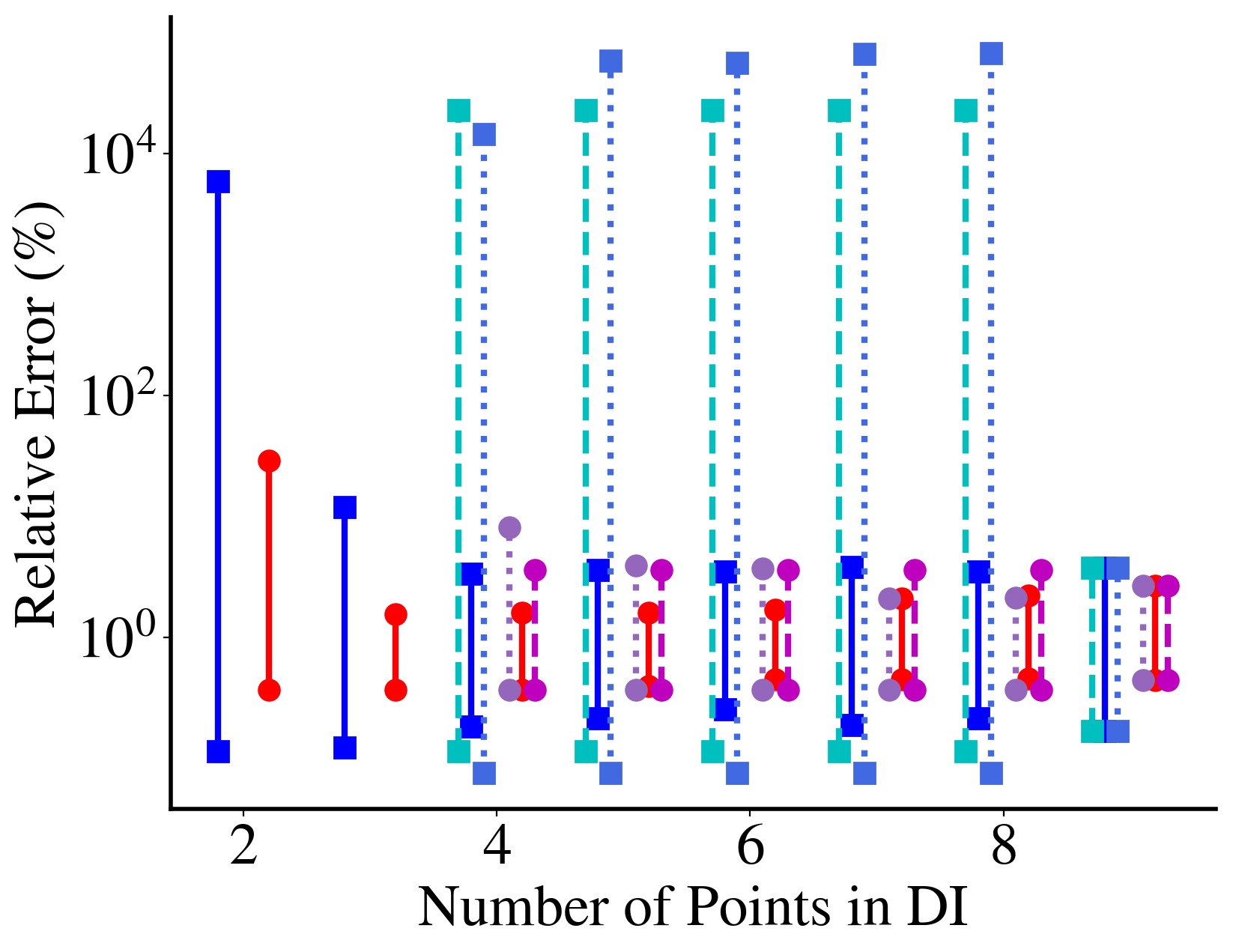}
        \caption{Error for degree 1 dynamical system in DI.}
    \end{subfigure}
    \begin{subfigure}[b]{.45\linewidth}
        \centering
        \includegraphics[width = \linewidth]{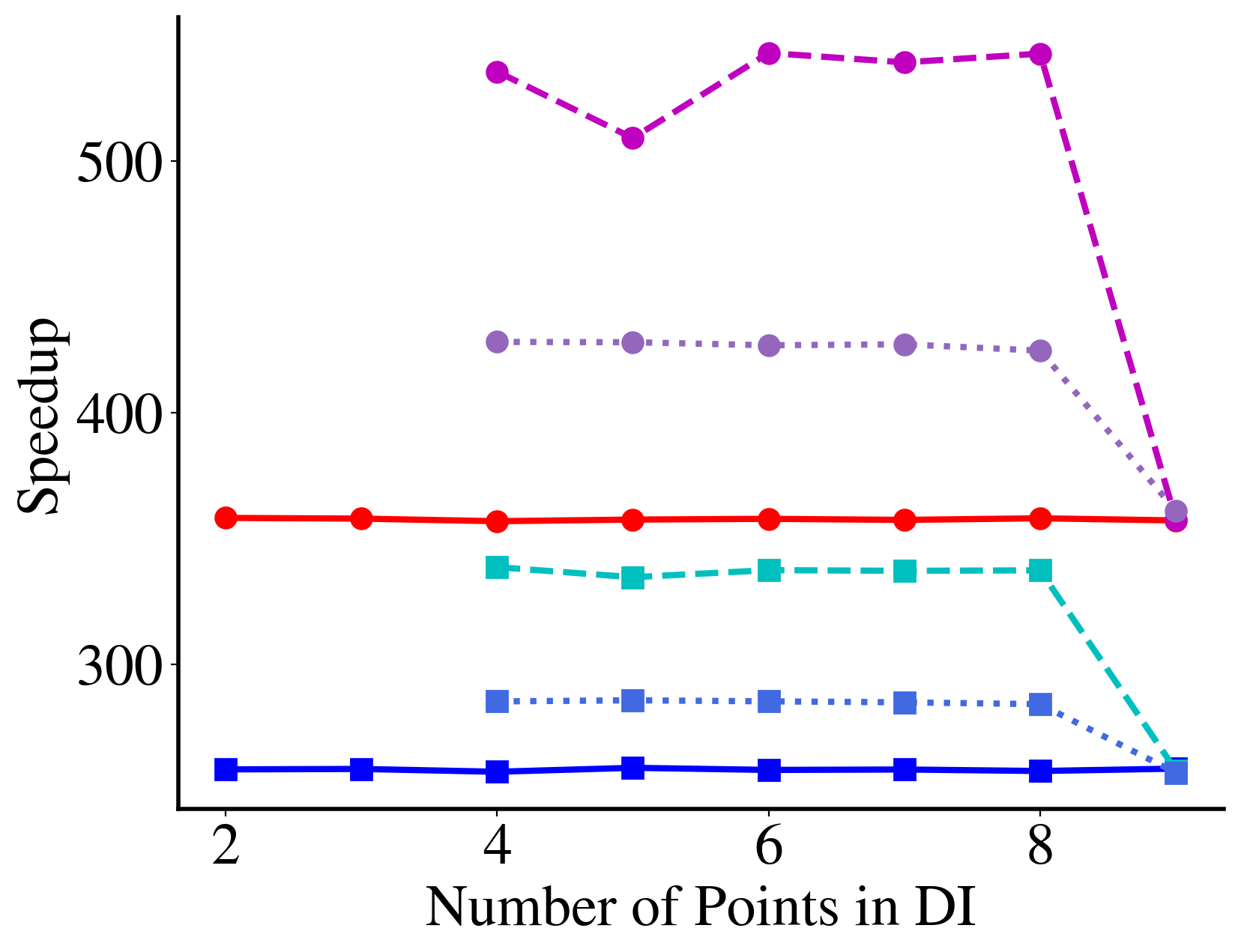}
        \caption{Speedup for degree 0 dynamical system in DI.}
    \end{subfigure}
    \begin{subfigure}[b]{.45\linewidth}
        \centering
        \includegraphics[width = \linewidth]{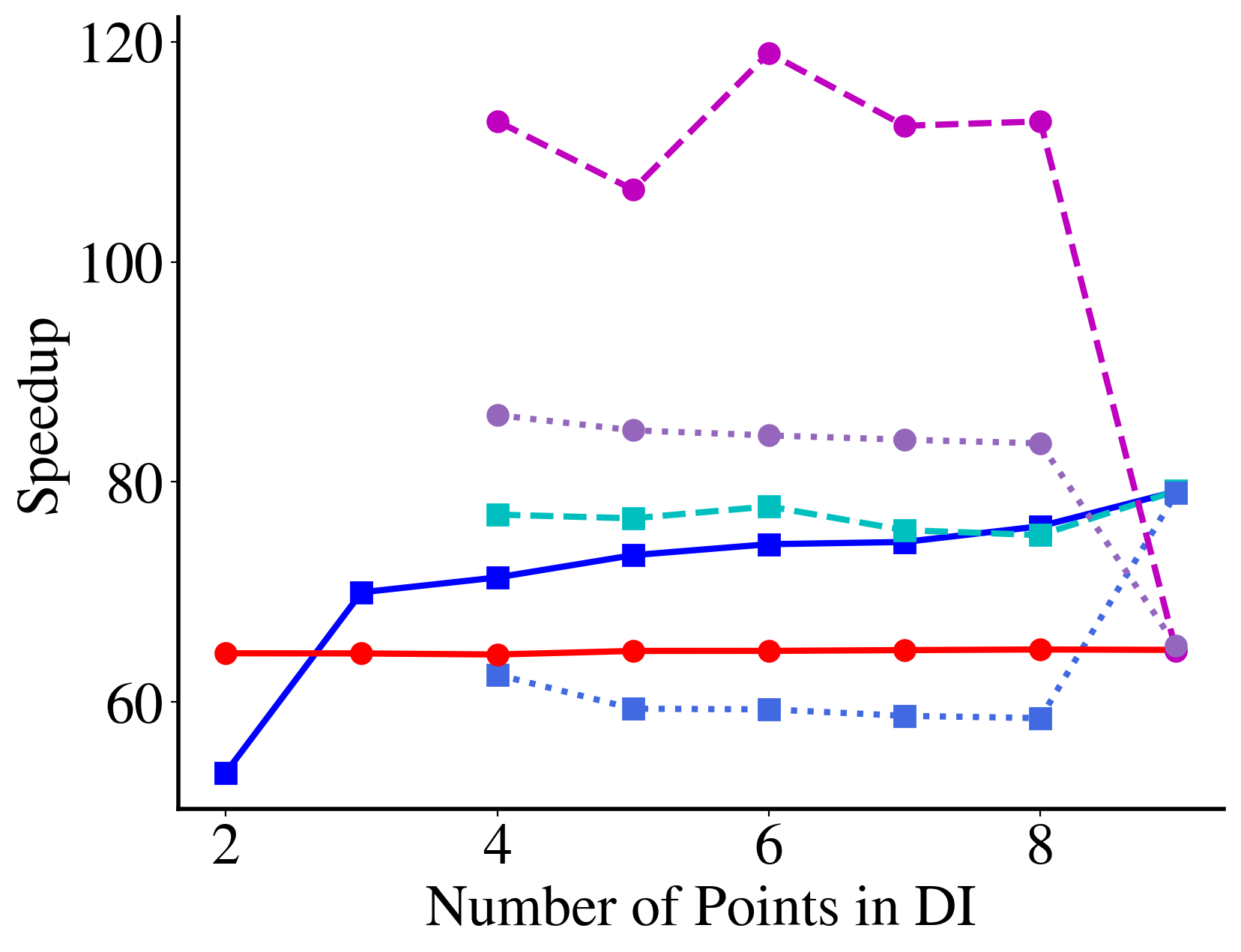}
        \caption{Speedup for degree 1 dynamical system in DI.}
    \end{subfigure}
\caption{Performance comparison among various models and the degree of polynomial in the DI. Nine training points are used in total. Both constant dynamical systems (left) and linear dynamical systems (right) are used in DI. The latent space dimension of four is used for LaSDI-NM models and the latent space dimension of five is used for LaSDI-LS models.  }
\label{fig: 1db compare}
\end{figure}

For more extensive analysis, we report accuracy (i.e., the relative error range) and speedup for various LaSDI models in Figure~\ref{fig: 1db compare}. We train LaSDI-NM and LaSDI-LS models with nine training points. Then both the region-based and interpolation-based local LaSDIs are constructed for various $\numDI$s, labeled as LaSDI-NM, interpolated LaSDI-NM with RBF, interpolated LaSDI-NM with B-Spline, LaSDI-LS, interpolated LaSDI-LS with RBF, and interpolated LaSDI-LS with B-Spline. Furthermore, we compare both Degree 0 and 1 dynamical system in DIs. There are several things to note. \textbf{First}, the relative ranges for Degree 1 dynamical systems tend to be bigger than the ones for Degree 0 dynamical systems. This makes sense because a lot more variations are possible with a higher degree, which gives a higher chance to be both the best and worst model. As a matter of fact, the best model with the minimum relative error range is found with Degree 1 Dynamical system (i.e., the region-based LaSDI-LS with $\numDI = 3$). \textbf{Second}, the interpolated LaSDIs tend to give a larger relative error range than the region-based local LaSDIs. We speculate that this is because the training points were not optimally chosen for the interpolation. This fact will be further investigated in the follow-up paper. \textbf{Third}, a higher speed-up is achieved by Degree 0 dynamical system. This also makes sense because the solution process of the latent space dynamics with less degrees will be faster. 

\begin{figure}
    \centering
    \begin{subfigure}[b]{.45\linewidth}
        \centering
        \includegraphics[width = \linewidth]{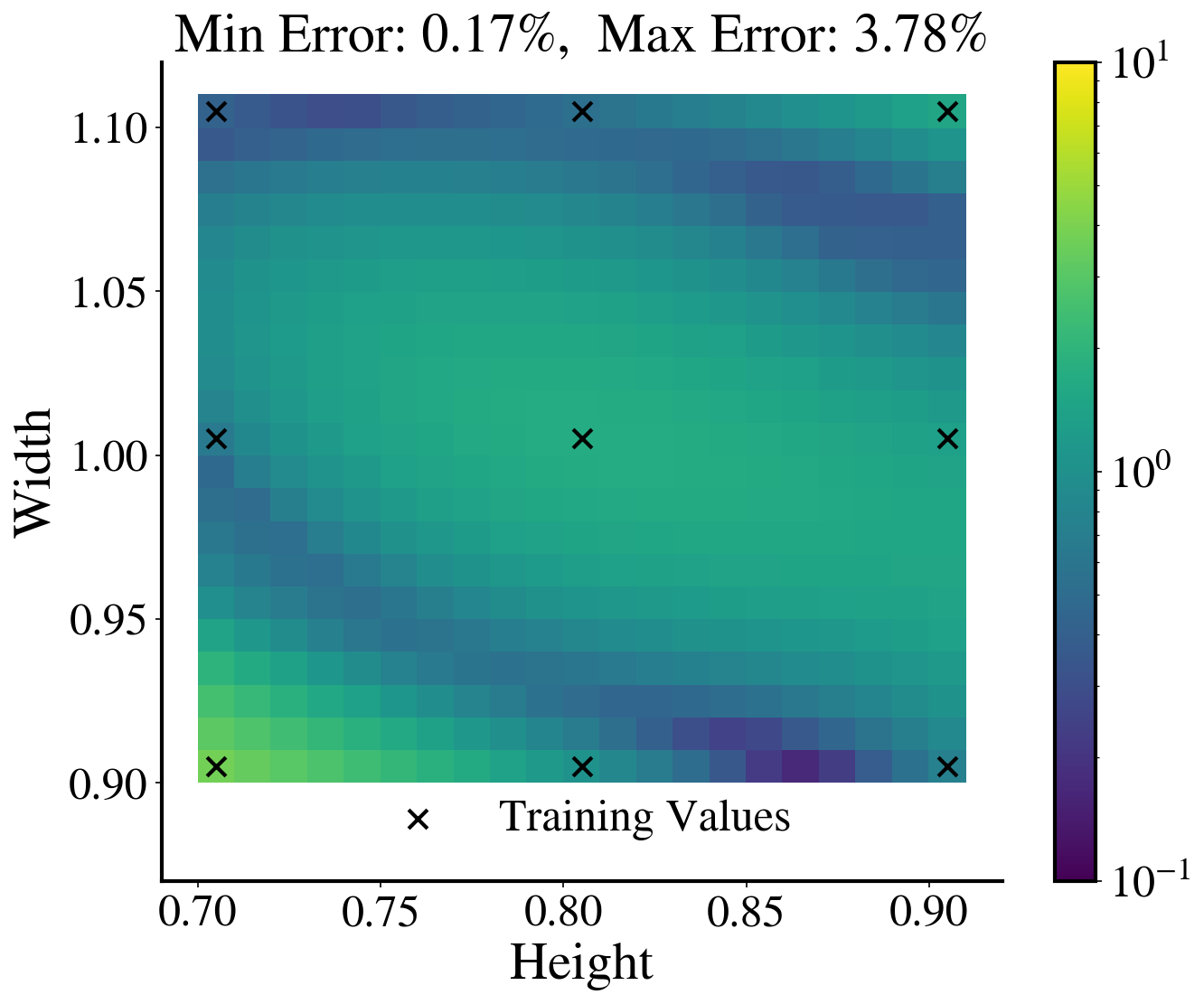}
        \caption{LaSDI-NM Global DI}
    \end{subfigure}
    \begin{subfigure}[b]{.45\linewidth}
        \centering
        \includegraphics[width = \linewidth]{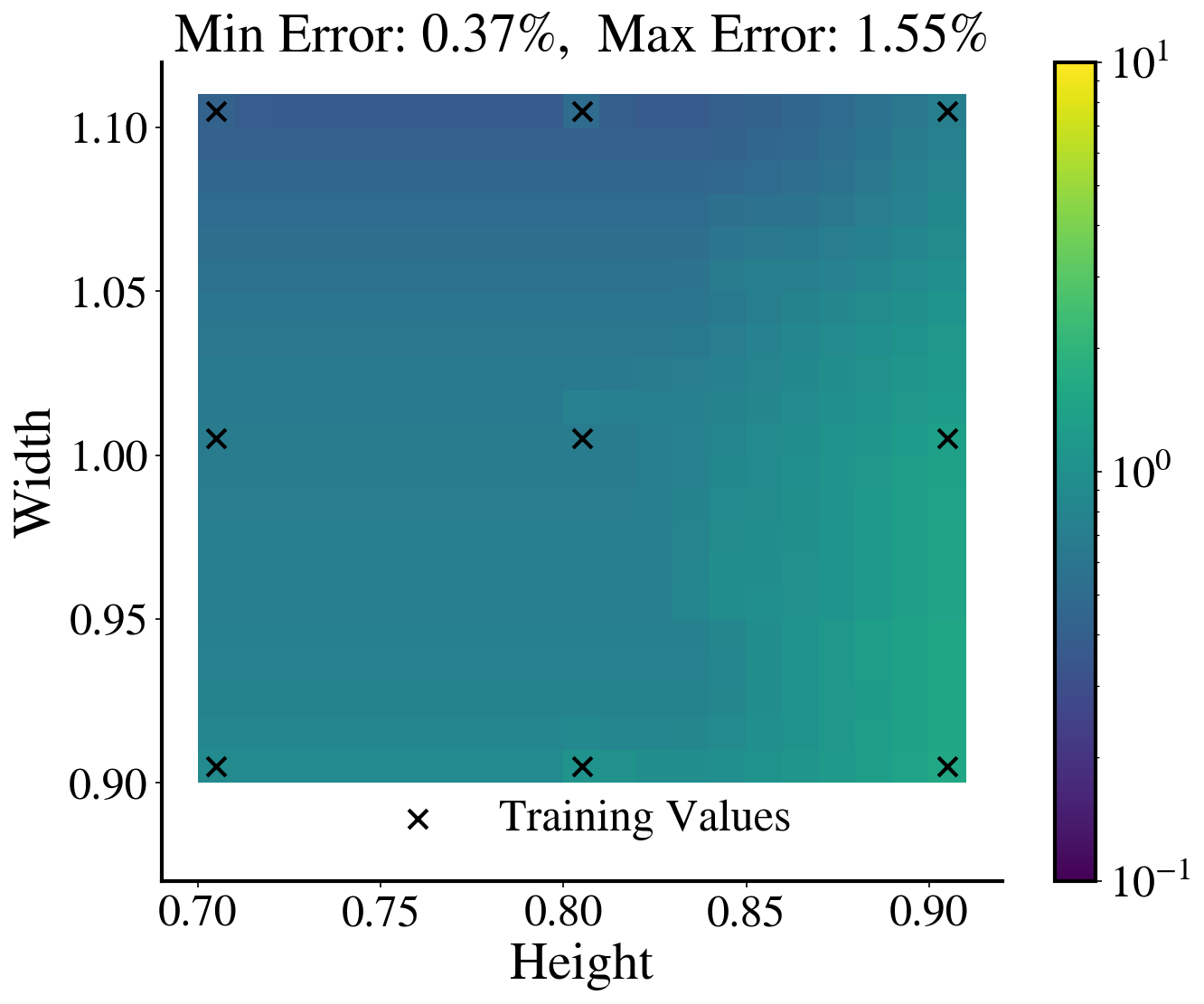}
        \caption{LaSDI-LS Local (3) DI}
    \end{subfigure}
    \caption{Heat maps from both LaSDI-NM with global DI and LaSDI-LS with local(3) DI that are identified from Figure~\ref{fig: 1db compare} in terms of accuracy for 1D Burgers problem. For both models, nine training points are used.}
    \label{fig:1db sum 9}
\end{figure}

\begin{figure}
    \centering
    \begin{subfigure}[b]{.45\linewidth}
        \centering
        \includegraphics[width = \linewidth]{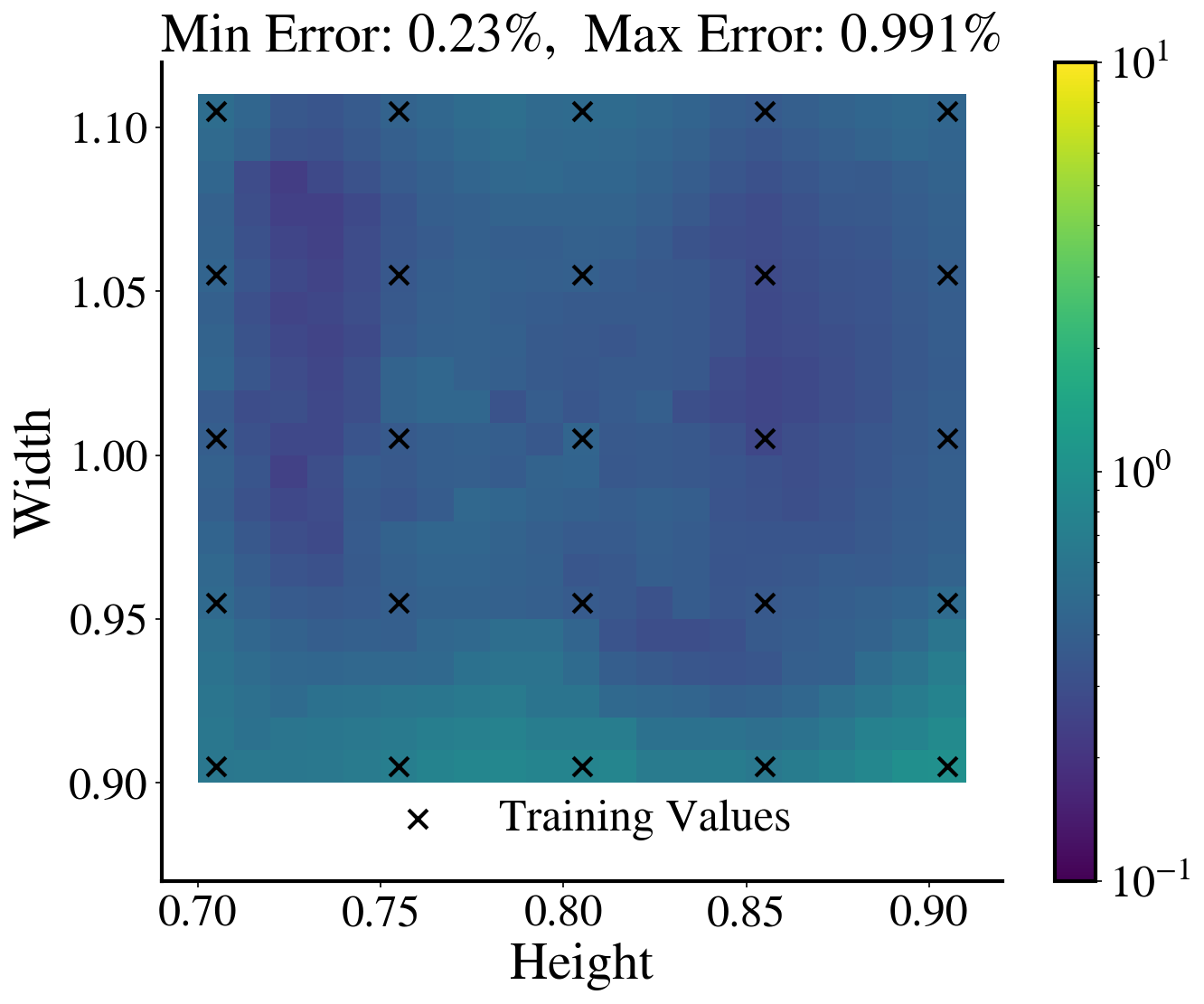}
        \caption{LaSDI-NM Local (16) DI}
    \end{subfigure}
    \begin{subfigure}[b]{.45\linewidth}
        \centering
        \includegraphics[width = \linewidth]{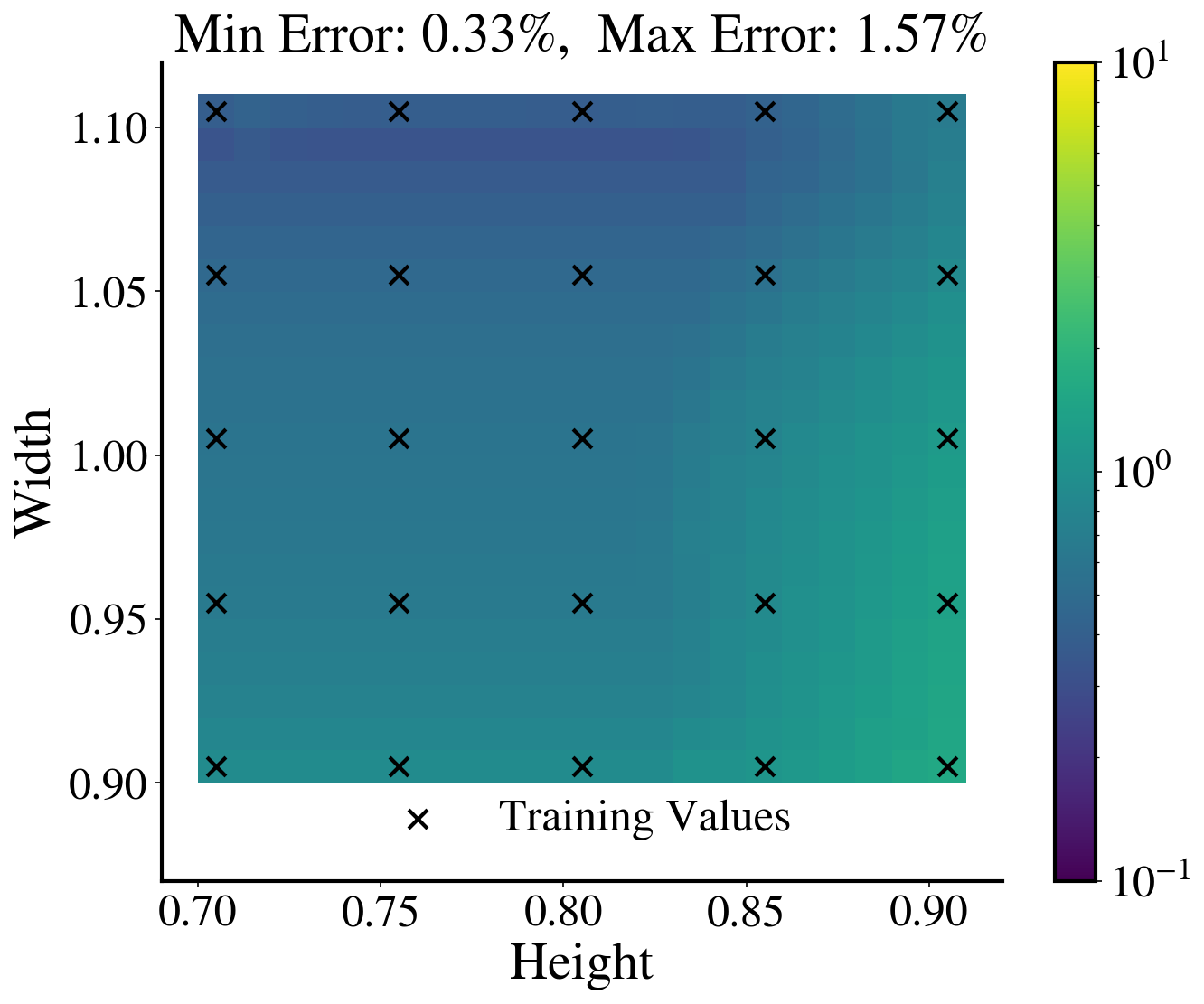}
        \caption{LaSDI-LS Local (5) DI}
    \end{subfigure}
    \begin{subfigure}[b]{\linewidth}
        \centering
        \includegraphics[width = \linewidth]{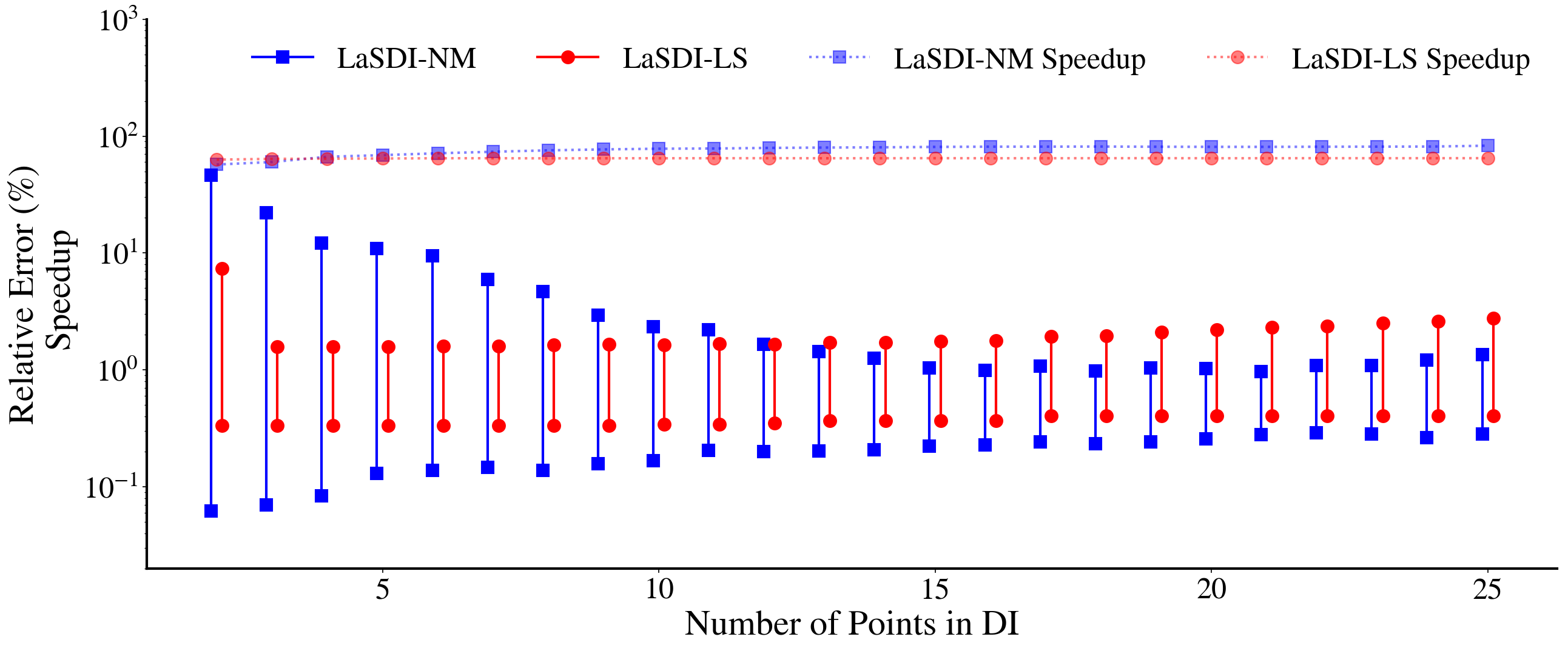}
        \caption{Relative Error Range and Speedup}
    \end{subfigure}
    \caption{Comparison between LaSDI-NM and LaSDI-LS on 1D Burgers problem with 25 training points. The region-based models are used for both LaSDI-NM and LaSDI-LS. For each case, the range of relative errors over all testing points in the parameter space are reported. We also show the heat maps of relative errors from both LaSDI-NM and LaSDI-LS for the corresponding best-case scenarios.}
    \label{fig:1db sum}
\end{figure}

In Figure~\ref{fig:1db sum 9}, we plot the heat maps of the relative error for the best LaSDI-LS and LaSDI-NM models that are identified from Figure~\ref{fig: 1db compare}, i.e., the LaSDI-NM Global DI and the region-based LaSDI-LS local DI with $\numDI=3$. The plot reveals that the LaSDI-LS local DI is better than the LaSDI-NM global DI model for this particular case.

Figure~\ref{fig:1db sum} compares the performance of LaSDI-LS and LaSDI-NM with 25 total training points to see the effect of increasing the number of training FOMs. First, the accuracy of LaSDI-NM has been improved when compared with the ones with 4 or 9 training points (see Figures~\ref{fig: 1db heat 4} and \ref{fig:1db sum 9}). For example, the maximum relative error of LaSDI-NM has decreased from $9.215\%$ to $3.78\%$ and to $0.991\%$ as the number of training points increase from 4 to 9 and to 25. This implying that more data is used, the better the accuracy of LaSDI-NM would be. Second, while LaSDI-NM gives a larger range of maximum relative errors when 9 training FOMs are used, it consistently gives the smallest relative error of all predictive LaSDIs. Finally, we notice no reduction in error in LaSDI-LS from 9 to 25 training values; this highlights the limitations of POD in advection-dominated problems.

\subsection{2D Burgers} \label{sec:2dburgersResults}
\begin{figure}
    \centering
    \begin{subfigure}[b]{.45\linewidth}
        \centering
        \includegraphics[width = \linewidth]{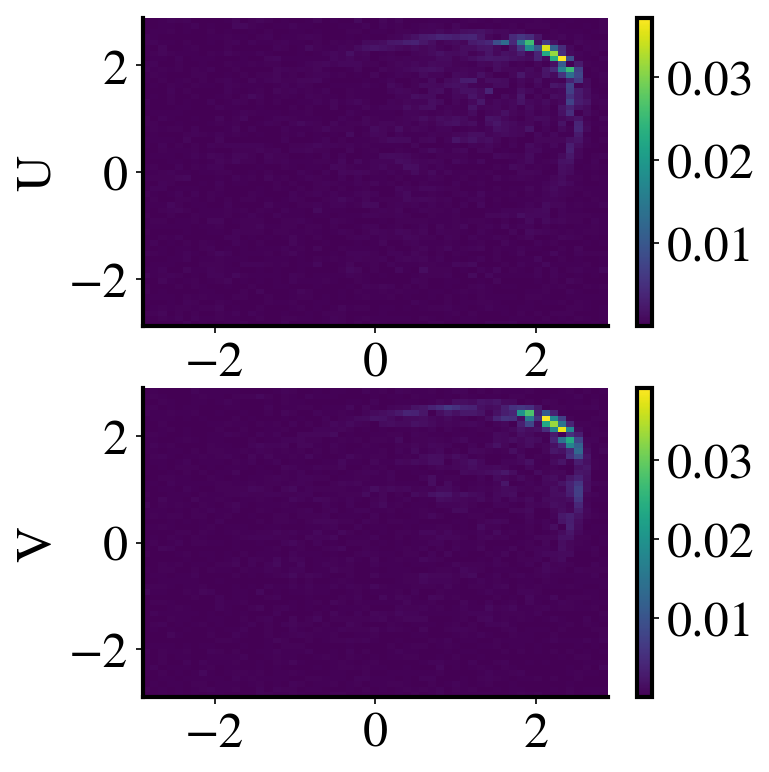}
        \caption{LaSDI-NM}
    \end{subfigure}
    \begin{subfigure}[b]{.45\linewidth}
        \centering
        \includegraphics[width = \linewidth]{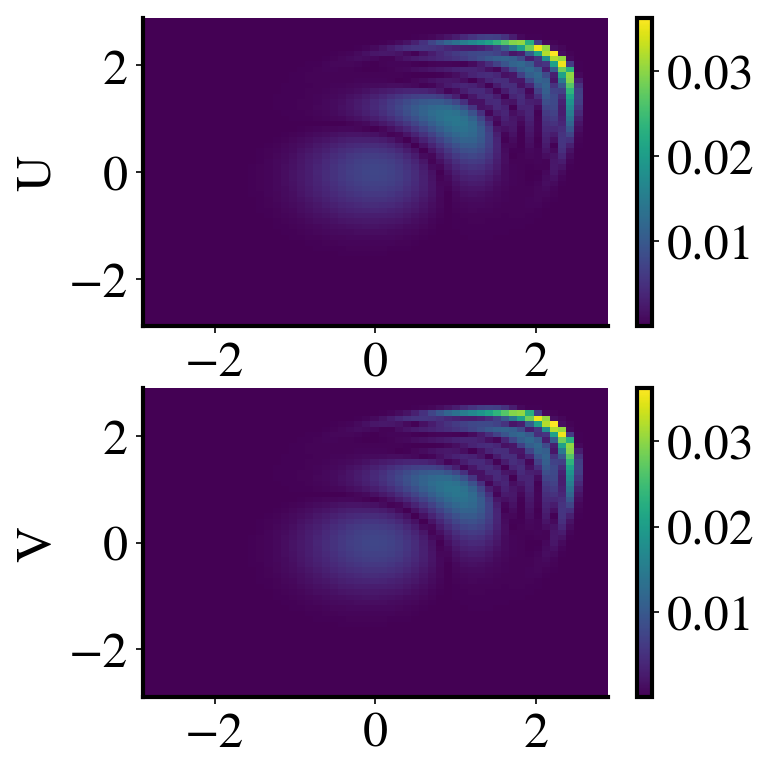}
        \caption{LaSDI-LS}
    \end{subfigure}
    \caption{The relative error of LaSDI models at the final time-step for 2D Burgers problem with $a = 0.8$ and $w = 1.02$. A latent space dimension of three is used for LaSDI-NM and a latent space dimension of five is used for LaSDI-LS.}
    \label{fig:2D single}
\end{figure}

We expect that LaSDI-NM will perform better than LaSDI-LS for the 2D Burgers problem based on the singular value decay shown in Figure~\ref{fig: SV Decay}. Figure~\ref{fig:2D single} illustrates the accuracy level achieved by LaSDI models for both the $u$ and $v$ components of the velocity as well as the relative error at the last time step. The figure  shows the clear advantage of LaSDI-NM over LaSDI-LS over the advection-dominated problem, as expected. Of course, the projection error of the linear compression can be improved by increasing the latent space dimension. However, the larger latent space dimension implies the more terms to be identified in the dynamics identification step. 

\begin{figure}
    \centering
    \includegraphics[width = \linewidth]{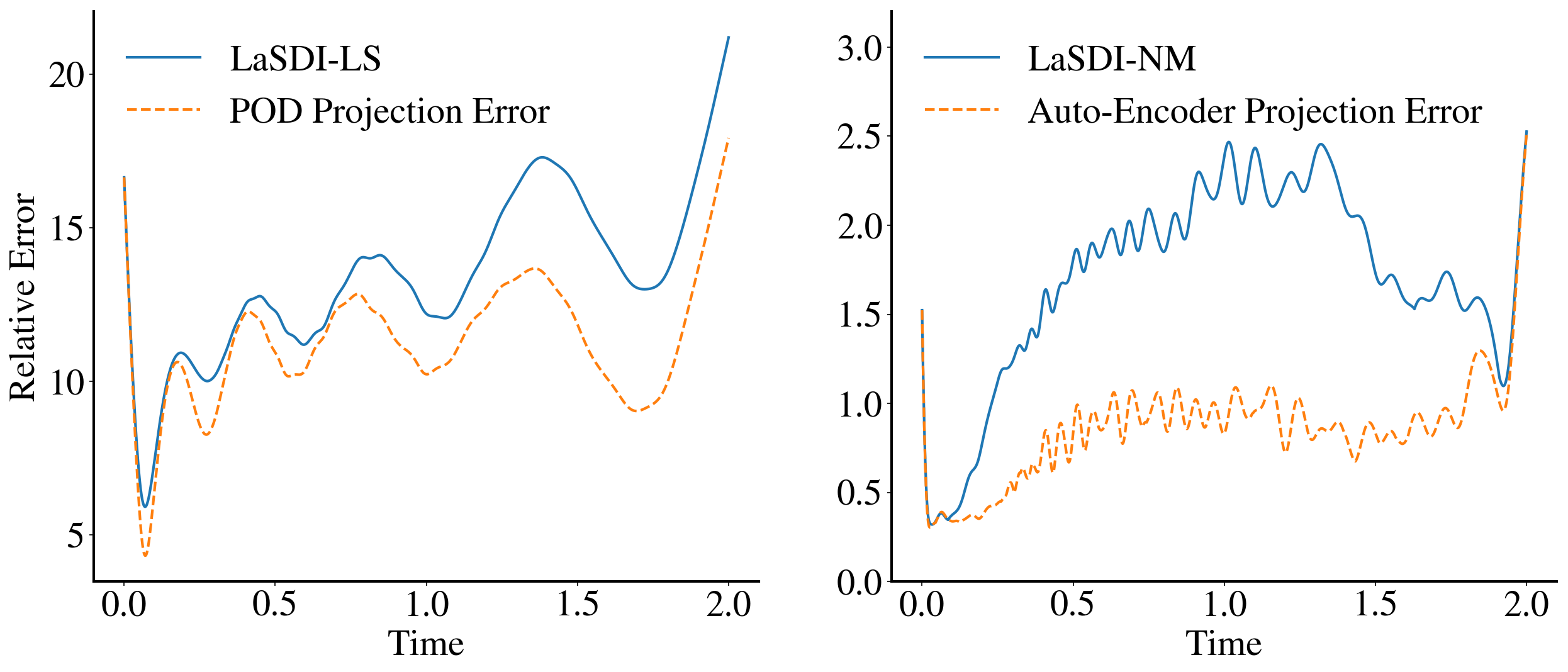}
    \caption{The illustration of relative errors of LaSDI models for the 2D Burgers simulation shown in Figure~\ref{fig:2D single}. Note the relative errors are bounded below by the projection error of each compression technique, i.e., proper orthogonal decomposition and auto-encoder. The projection errors are computed by compressing and de-compressing the corresponding full order model solutions either by POD basis or auto-encoder.}
    \label{fig:2d ae error}
\end{figure}

Figure \ref{fig:2d ae error}  illustrates the relative error of a LaSDI model is bounded below by the projection error of the compression technique. 

\begin{figure}
    \centering
    \begin{subfigure}[b]{.45\linewidth}
        \centering
        \includegraphics[width = \linewidth]{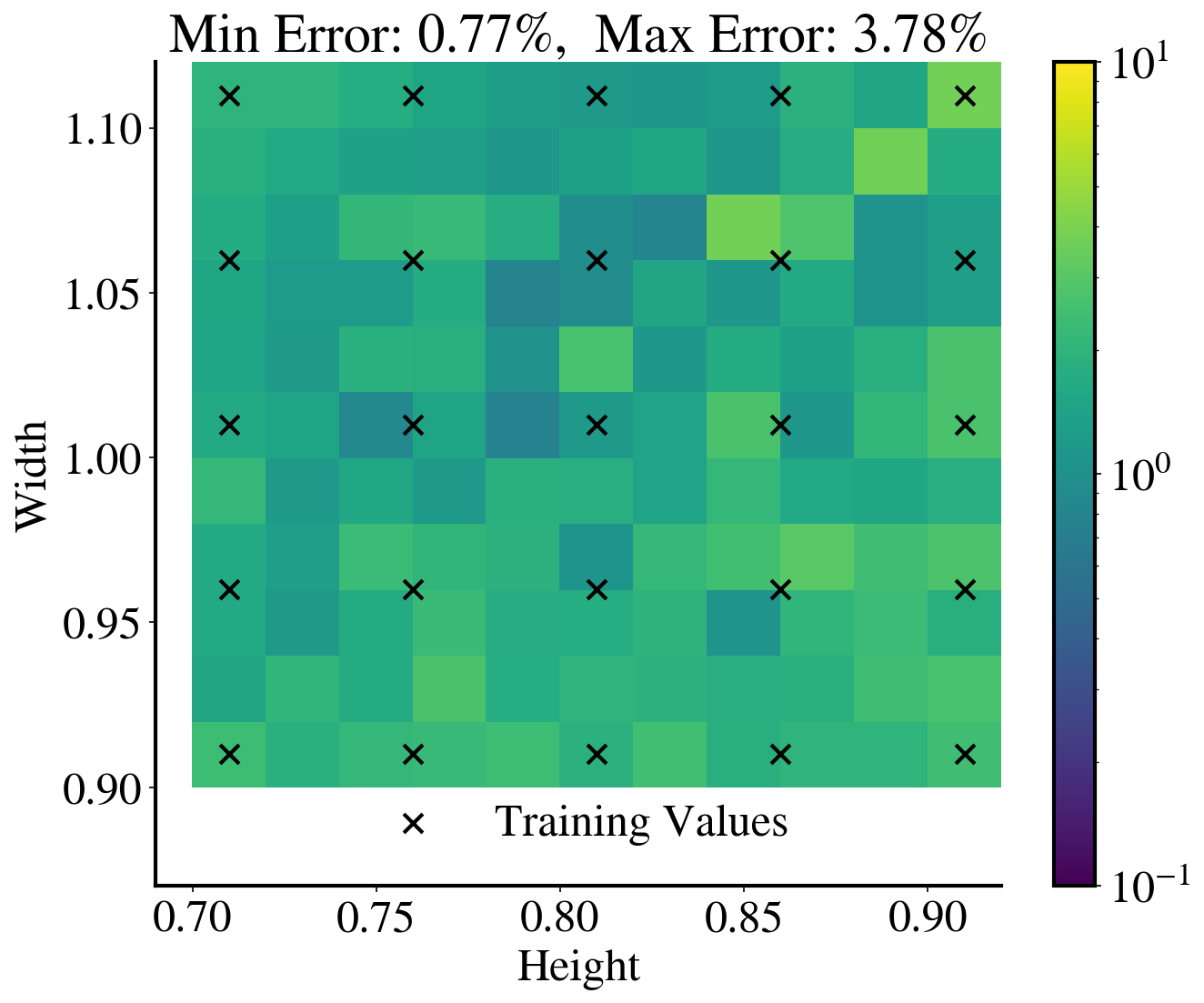}
        \caption{LaSDI-NM Local (6) DI}
    \end{subfigure}
    \begin{subfigure}[b]{.45\linewidth}
        \centering
        \includegraphics[width = \linewidth]{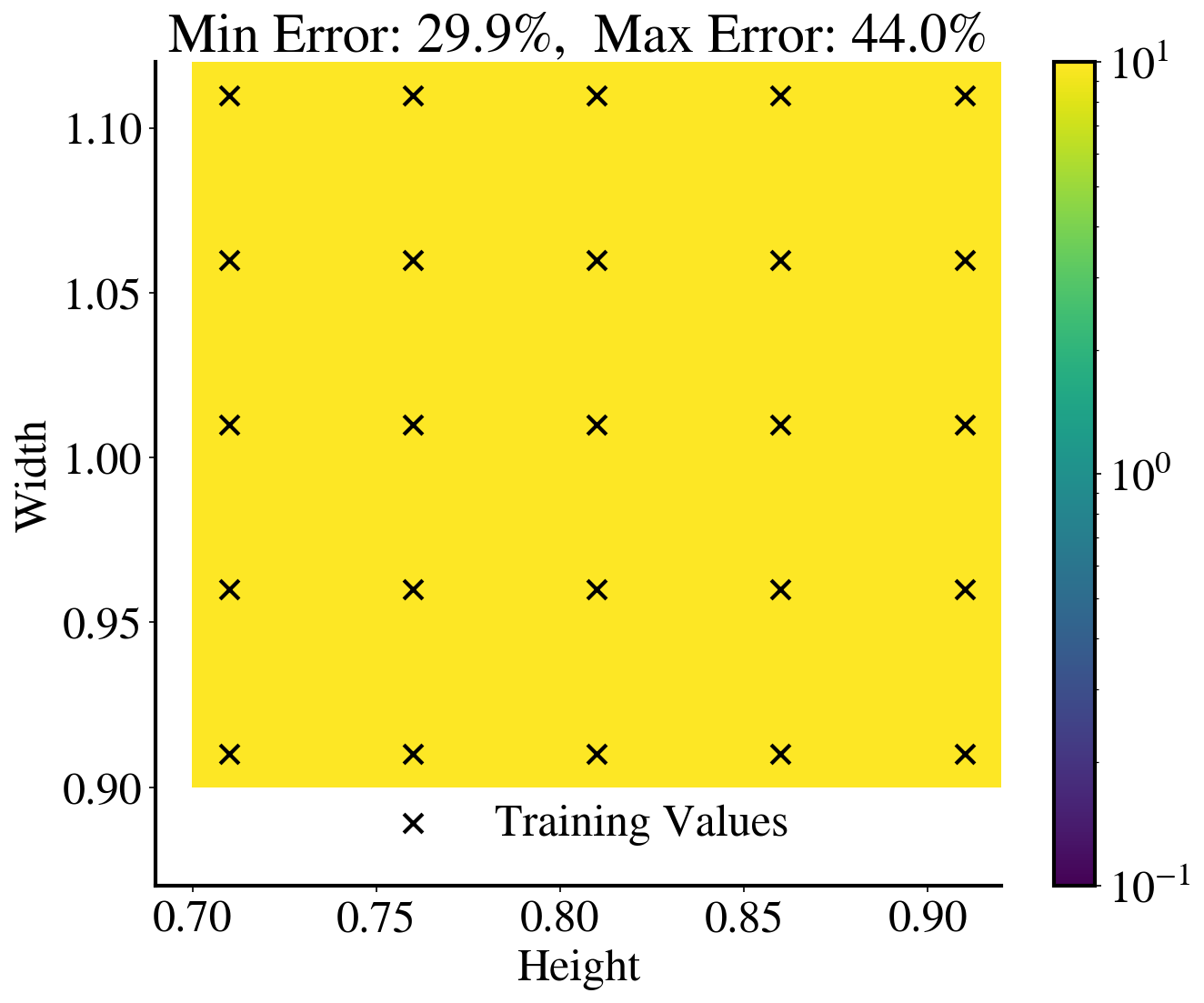}
        \caption{LaSDI-LS Global DI}
    \end{subfigure}
    \begin{subfigure}[b]{\linewidth}
        \centering
        \includegraphics[width = \linewidth]{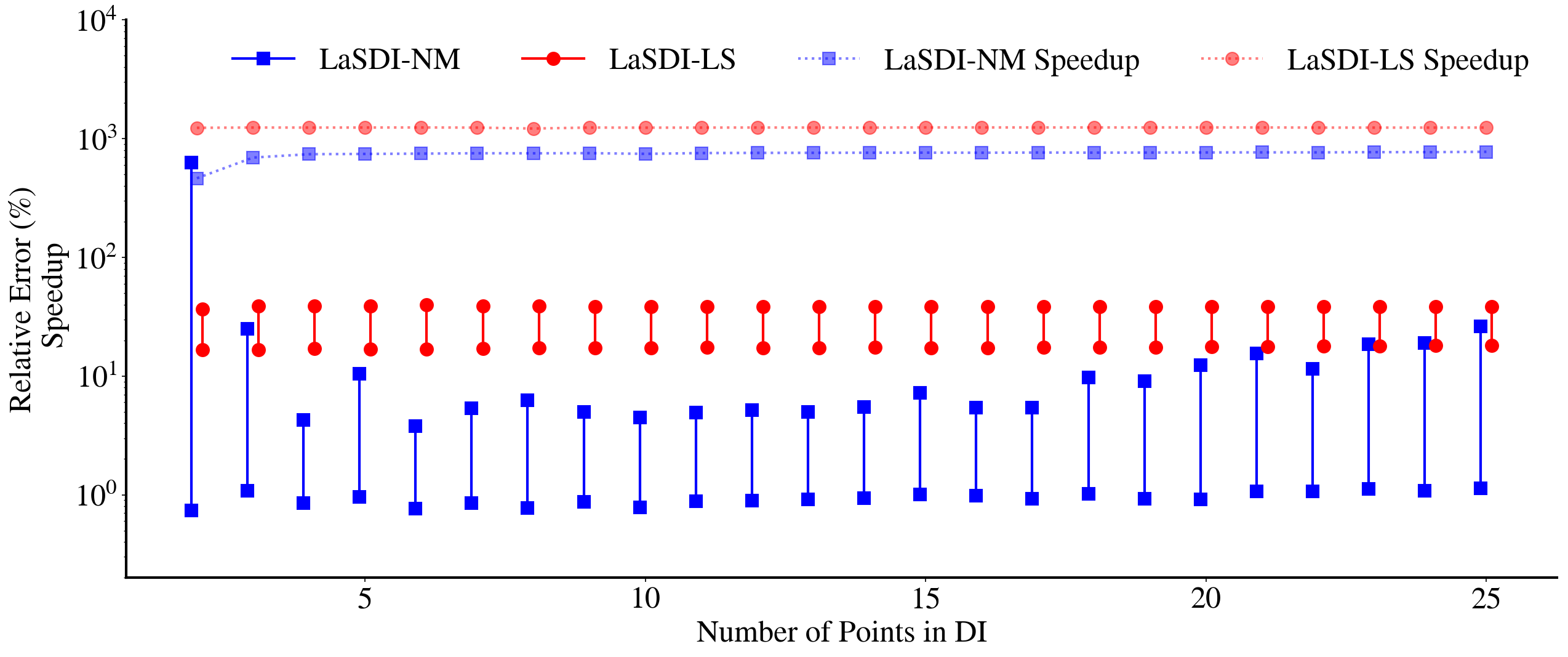}
        \caption{Relative Error Range and Speedup}
    \end{subfigure}
    \caption{Comparison between LaSDI-NM and LaSDI-LS on 2D Burgers problem with 25 training points. We include two heat maps of the relative error of the LaSDIs for the prescribed parameter values. LaSDI-NM uses a latent-space dimension of three and a cubic dynamical system in the DI. LaSDI-LS uses a five dimensional latent-space with a linear dynamical system. }
    \label{fig:2d GL compare}
\end{figure}

The comparison between LaSDI-LS and LaSDI-NM becomes even more vivid in Figure \ref{fig:2d GL compare}, which shows the performance of LaSDI models with 25 training points. For this example, we note that the region-based local DI  outperforms global DI for LaSDI-NM. This implies that the latent-space dynamics are heavily localized within the parameter space. That is, for example, the dynamical system that best approximates the case of $a = 0.6$, $w = 0.8$ is not a good model for the case of $a = 0.9$, $w = 1.1$.

For this particular problem, a speed-up around 800 is achieved by LaSDI-NM. We use a cubic dynamical system, with cross-terms, to improve the expressivity and therefore accurately represent the three dimensional latent space generated in LaSDI-NM. However, the higher-degree dynamical systems in DI can be numerically unstable. Furthermore, the lower degree approximation can lead to a higher speed-up. Thus, the balance among expressivity, numerical instability, and speed-up must be sought when LaSDI models are applied.

\subsection{Time-Dependent Heat Conduction}

\begin{figure}
    \centering
    \begin{subfigure}[b]{.435\linewidth}
        \centering
        \includegraphics[width = \linewidth]{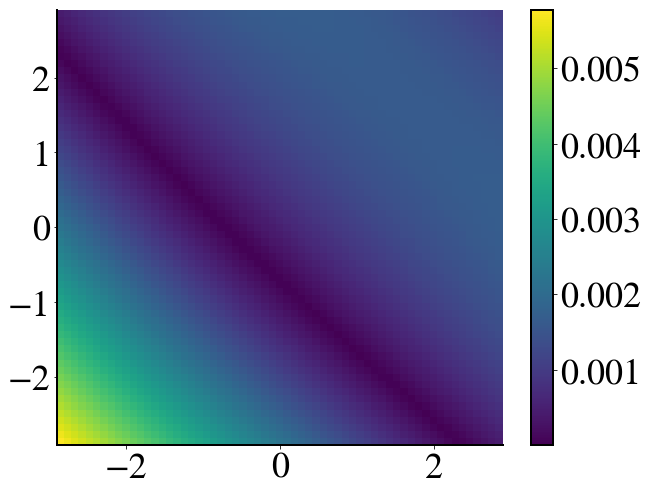}
        \caption{LaSDI-NM}
    \end{subfigure}
    \begin{subfigure}[b]{.465\linewidth}
        \centering
        \includegraphics[width = \linewidth]{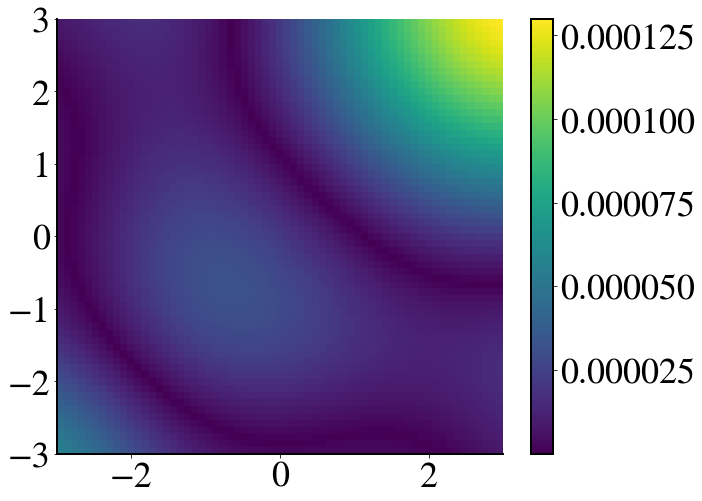}
        \caption{LaSDI-LS}
    \end{subfigure}
    \caption{Relative error at the final time step for $w=1.0$ and $a=2.0$ in the heat conduction problem generated using both LaSDI-NM and LaSDI-LS. LaSDI-NM uses a three dimensional latent-space whereas LaSDI-LS uses a four dimensional latent-space. In both cases, a linear dynamical system is used in DI.}
    \label{fig:Ex16 single}
\end{figure}

Based on the singular value decay in Figure~\ref{fig: SV Decay}, we expect that the LaSDI-LS will perform better than the LaSDI-NM. Indeed, Figure~\ref{fig:Ex16 single} shows that the LaSDI-LS achieves much smaller relative error than the LaSDI-NM for $w=1.0$ and $a=2.0$.

\begin{figure}
\centering
    \includegraphics[width = .5\linewidth]{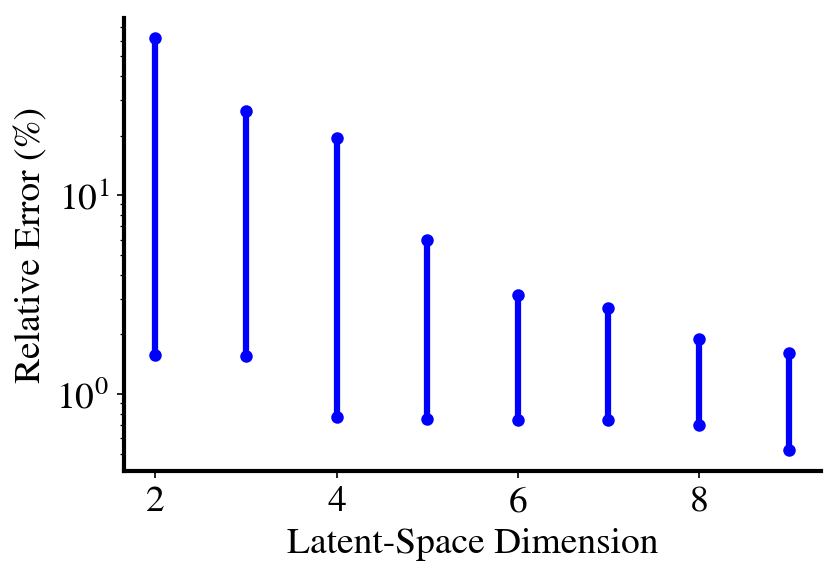}
    \caption{Relative error ranges for LaSDI-LS applied to the heat conduction problem with various latent space dimensions. The training parameters, $\trainingSet$,  are used to build LaSDI-LS models and testing parameters, $\testingSet$, are used to compute the relative error range, as prescribed in Table~\ref{tab:training}. In each case, we use Global DI with a linear dynamical system.}
    \label{fig:Ex16 sv error}
\end{figure}

As discussed in Section~\ref{sec:POD}, increasing the dimension of the the latent space necessarily increases the amount of information retained in the data-compression. To quantify this, we track the first $n_s$ singular values as a proportion of the sum of the singular values:
\begin{equation}\label{eq:svmass}
    m_{\text{sv}} = \sum_{i = 1}^{n_s} \sigma_i/\sum_{i=1}^{N_s}\sigma_i,
\end{equation}
where $m_{\text{sv}}$ serves as an indicator for the projection error of the POD-based reduced space as the projection error decreases as $m_{\text{sv}}$ increases. As seen in Figures~\ref{fig: 1db errors} and \ref{fig:2d ae error}, the projection error of the POD also serves as the lower bound for the LaSDI-LS. Therefore, it would be nice if the accuracy of the LaSDI-LS can follow the trend of the projection error of the POD-based compression. Figure~\ref{fig:Ex16 sv error} indeed illustrates the relative error range decreases as the latent space dimension increases. This feature implies that the accuracy of LaSDI-LS follows the trend of the projection error for this particular problem.

The decrease in the projection error implies that the quality improvement of the latent space accomplished by the POD. However, the overall accuracy of LaSDI models does not depend only on the quality of the latent space, but also on the quality of the DI model. Therefore, the decay of the relative error range in Figure~\ref{fig:Ex16 sv error} is possible because the DI models are in good quality as well. 

\begin{figure}
    \centering
    \begin{subfigure}[b]{.45\linewidth}
        \centering
        \includegraphics[width = \linewidth]{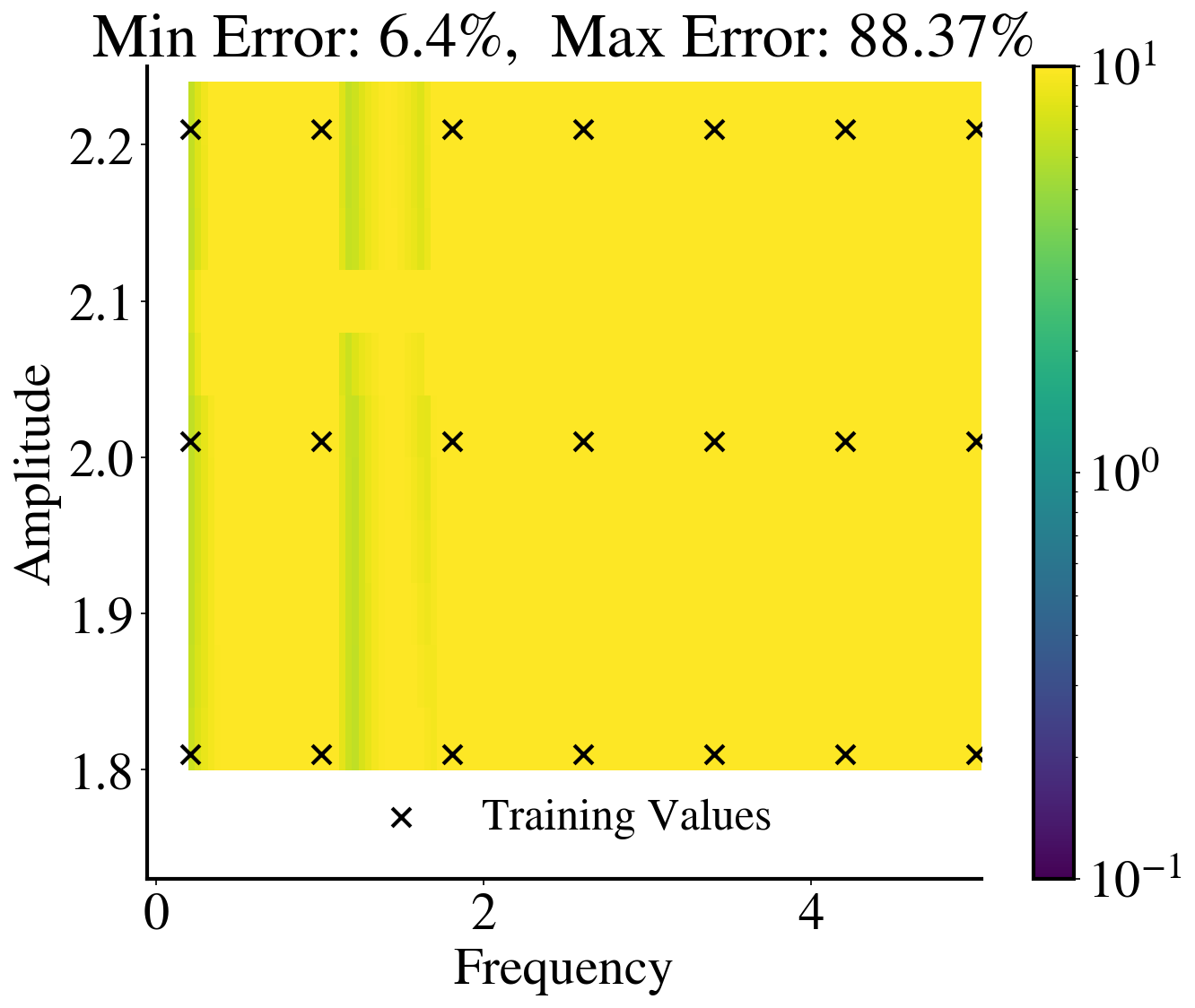}
        \caption{LaSDI-NM Global DI}
    \end{subfigure}
    \begin{subfigure}[b]{.45\linewidth}
        \centering
        \includegraphics[width = \linewidth]{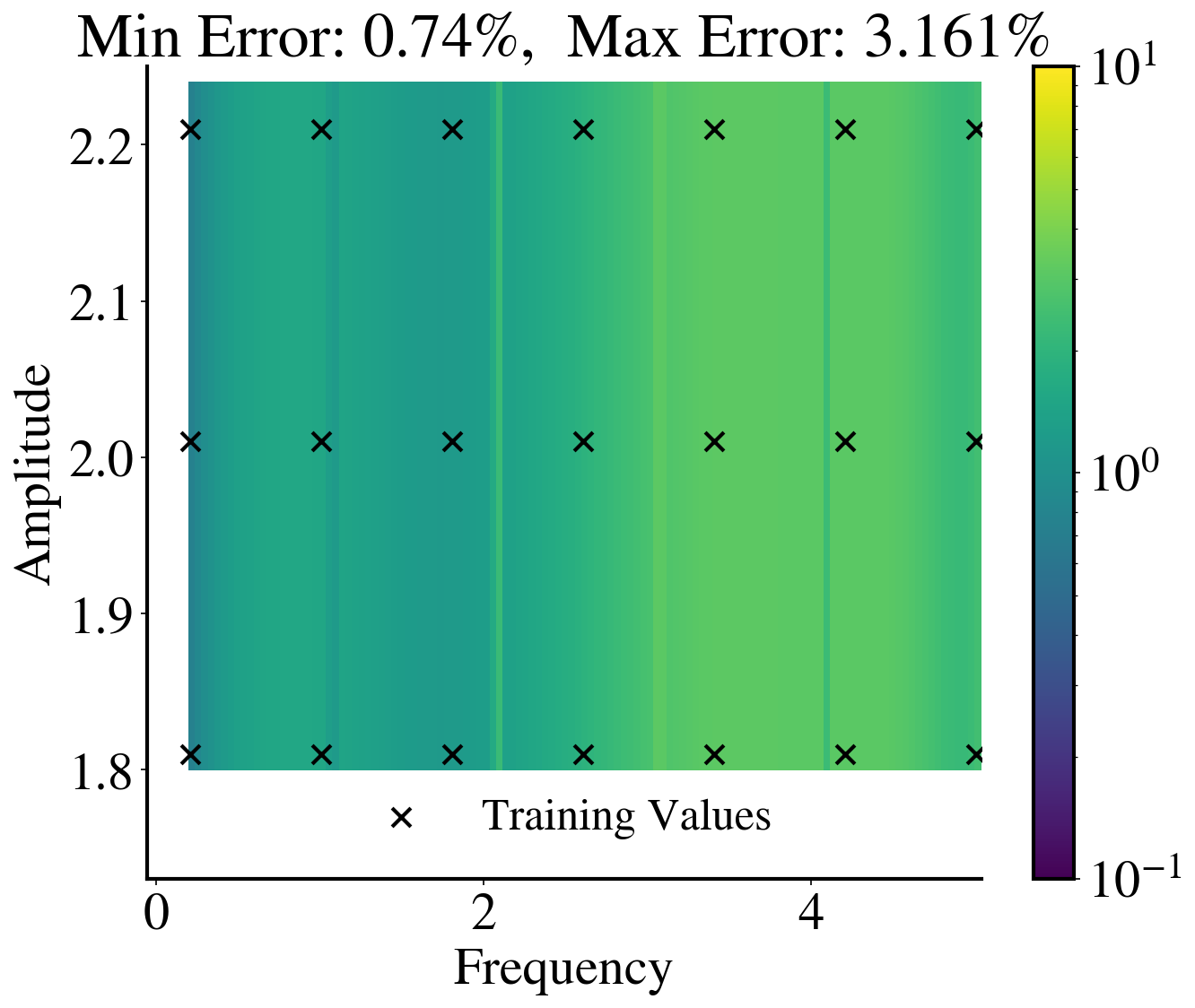}
        \caption{LaSDI-LS Global DI}
    \end{subfigure}
    \begin{subfigure}[b]{\linewidth}
        \centering
        \includegraphics[width = \linewidth]{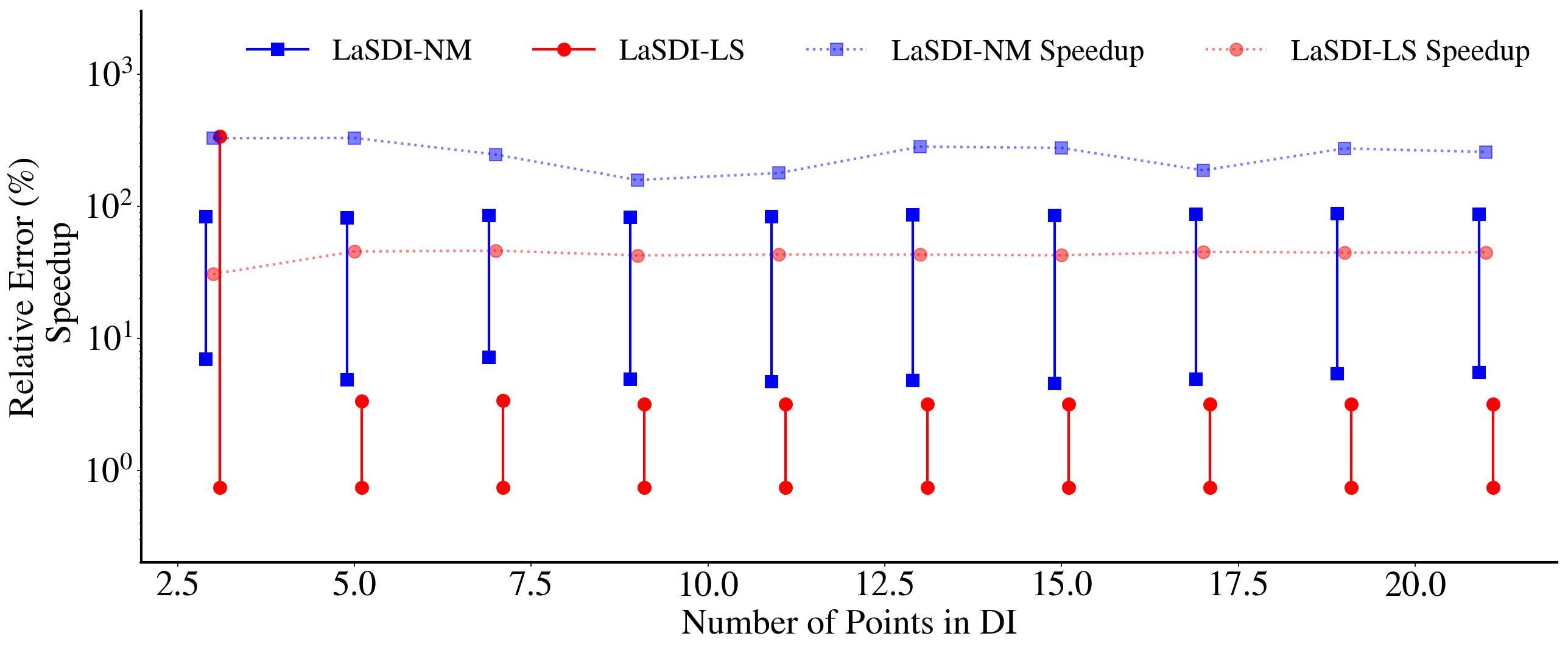}
        \caption{Relative Error Range and Speedup}
    \end{subfigure}
    \caption{Comparison between LaSDI-NM and LaSDI-LS for the Heat Conduction problem with 21 training values. We include two heat maps of the maximum relative error of the LaSDI models for the prescribed parameter values when Global DI is implemented. For LaSDI-NM and LaSDI-LS, we use a constant and linear dynamical system, respectively, during the DI process.}
    \label{fig:ex16 GL compare}
\end{figure}

Accurate LaSDI-LS requires fast singular value decay and as we present below, simple latent space dynamics are required to maintain this accuracy (i.e., Figure~\ref{fig:ex16 GL compare} demonstrates that 0 degree dynamical systems in DI is enough to produce a good accuracy for LaSDI-LS). Although the computational cost for the integration of the ODE and the matrix multiplication to decompress the data increases as we increase $n_s$, the decrease in LaSDI-LS relative error at the expense of decreased speedup times clearly gives benefits to applications that require a high accuracy.  

When expanding across the full parameter space, LaSDI-LS continues to produce lower LaSDI-LS errors than LaSDI-NM, as shown in Figure \ref{fig:ex16 GL compare}. The LaSDI-NM does not produce accurate results leading to the maximum $88\%$ error of LaSDI-NM, while LaSDI-LS achieves the maximum $3.16\%$ error. Due to the simplicity of both the linear and nonlinear latent spaces, we do see consistent results between local and global DI. As in Section \ref{sec:2dburgersResults}, simpler dynamics (i.e., 0-degree dynamical system) was enough to approximate the nonlinear latent space dynamics accurately, also leading to a speed-up of much higher than $100$x. 

While further neural network architectures can be explored, we restrict ourselves to the shallow masked network, which gave good performance for the 2D Burgers problem in Section~\ref{sec:2dburgersResults}. Because we expect diffusion-dominated problems to be well-represented in a latent spaces that is constructed by a linear compression, as suggested in Figure~\ref{fig: SV Decay}, exploring various neural network structures for this specific problem becomes counter-productive. Furthermore, nonlinear compression techniques, such as auto-encoders, require more computationally expensive training phase than linear compression techniques, such as POD. We have included the results of LaSDI-NM only to highlight the importance of data compression selection when implementing LaSDI.

\subsection{Radial Advection}\label{sec:radAdvection}

\begin{figure}
    \centering
    \begin{subfigure}[b]{.45\linewidth}
        \centering
        \includegraphics[width = \linewidth]{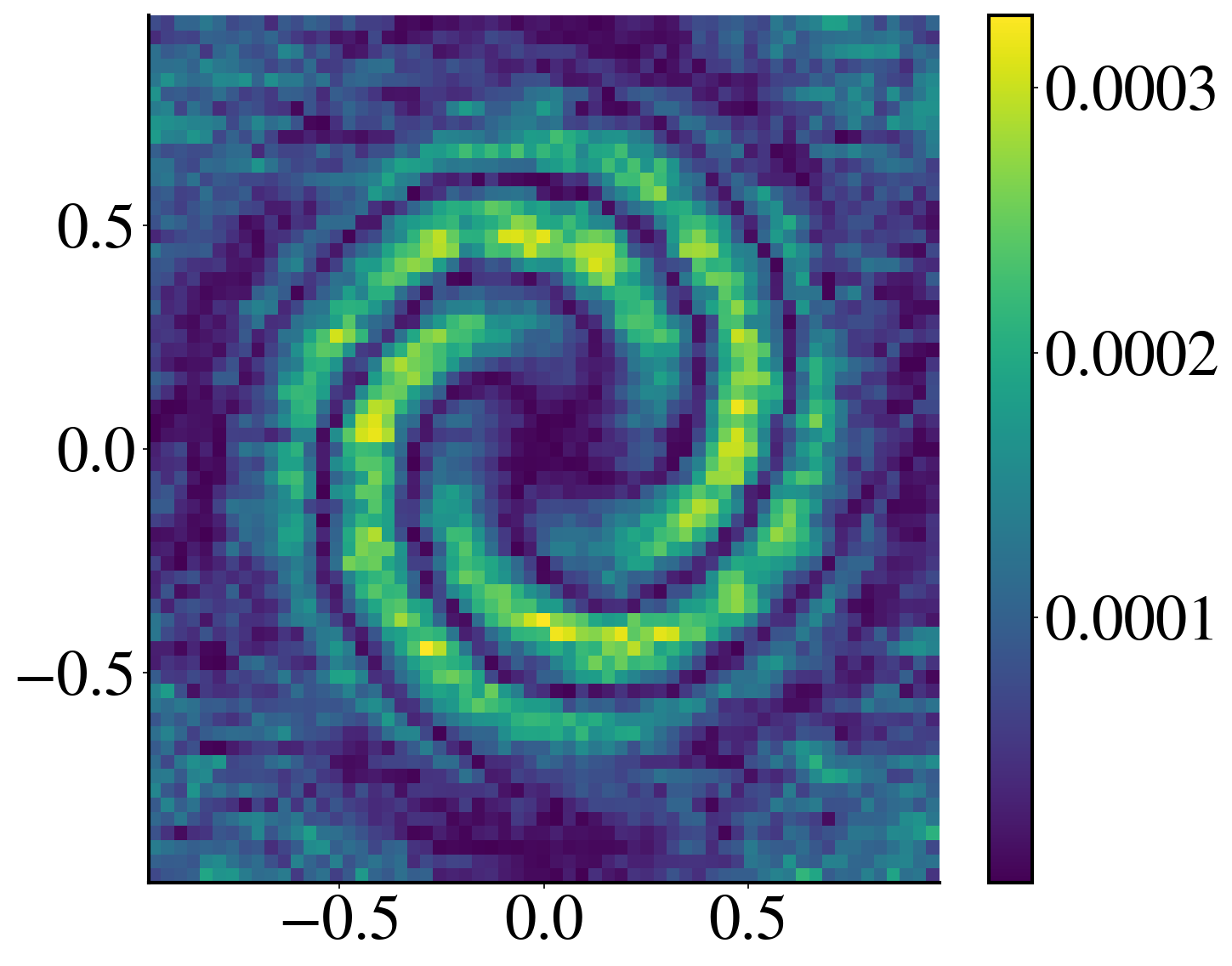}
        \caption{LaSDI-NM}
    \end{subfigure}
    \begin{subfigure}[b]{.45\linewidth}
        \centering
        \includegraphics[width = \linewidth]{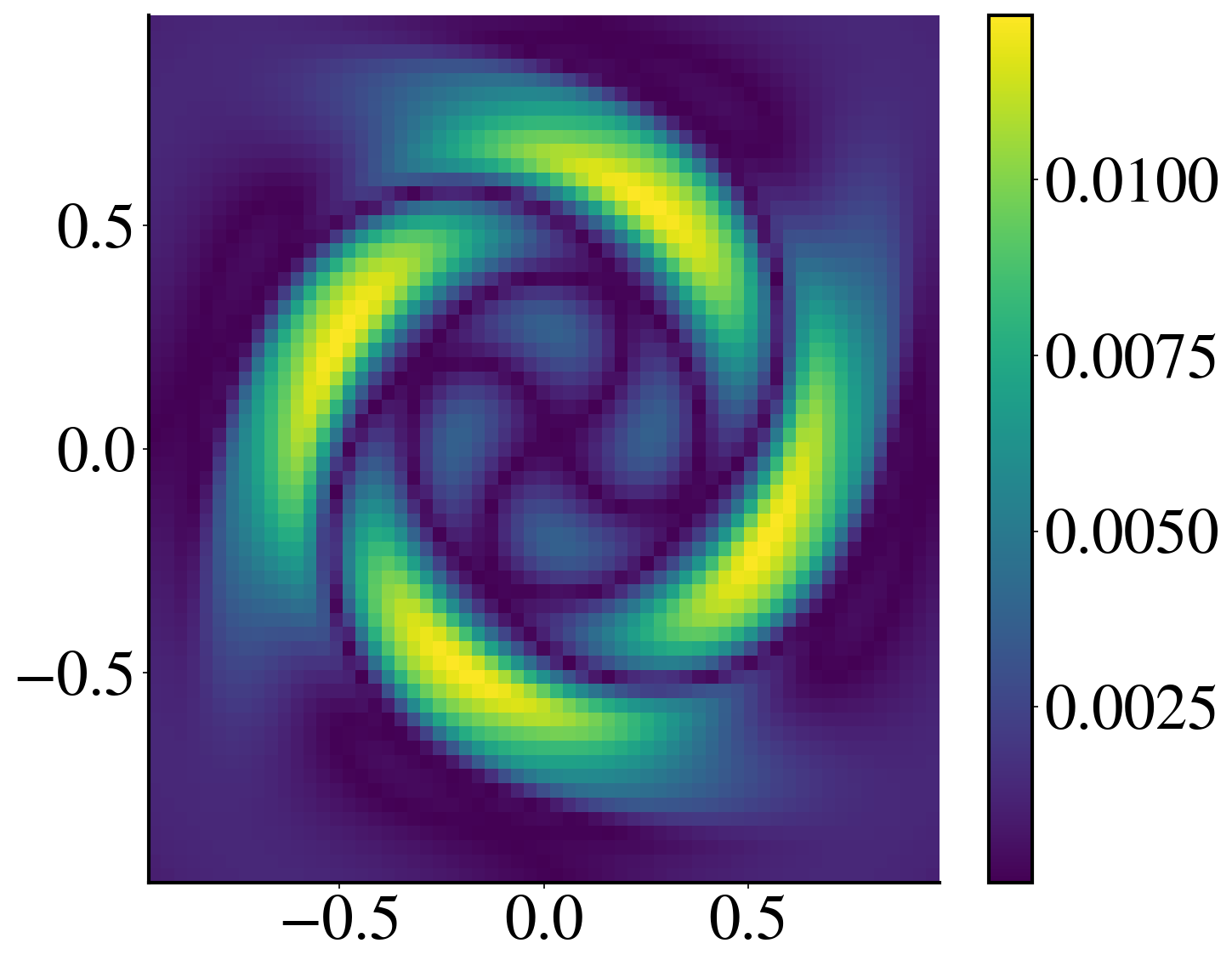}
        \caption{LaSDI-LS}
    \end{subfigure}
    \caption{Heat maps of the relative error for LaSDI-NM and LaSDI-LS on the radial advection problem for $\omega = 1.01$. For LaSDI-LS, we use a latent-space dimension of four with a linear dynamical system in DI. For LaSDI-NM, we prescribe a latent-space dimension of three with a quadratic dynamical system and no cross terms in DI. In both cases, global DI is used.}
    \label{fig:ex9 single}
\end{figure}

We show the numerical results for the radial advection problem. As in the previous sections, we first present results of LaSDI models for one specific test point to illustrate the viability of the method. Figure \ref{fig:ex9 single} presents the heat maps of relative errors produced by both LaSDI-LS and LaSDI-NM. As expected, LaSDI-NM performs significantly better than LaSDI-LS in terms of accuracy because the problem is advection-dominated and the singular value decay is slow as seen in Figure~\ref{fig: SV Decay}.

\begin{figure}
    \centering
    \begin{subfigure}[b]{\linewidth}
        \centering
        \includegraphics[width = \linewidth]{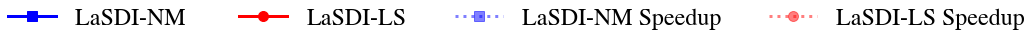}
    \end{subfigure}
    \begin{subfigure}[b]{.45\linewidth}
        \centering
        \includegraphics[width = \linewidth]{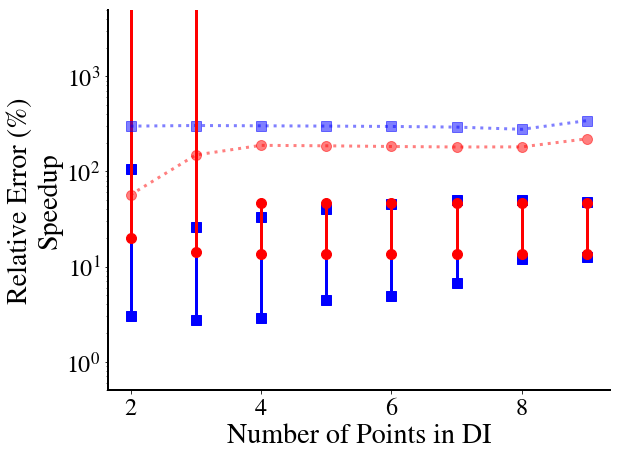}
        \captionsetup{justification=centering}
        \caption{$\trainingSet = \{0.60,0.70,...,1.40\}$\\ $\testingSet = \{0.60,0.61,...,1.40\}$}
    \end{subfigure}
        \begin{subfigure}[b]{.45\linewidth}
        \centering
        \captionsetup{justification=centering}
        \includegraphics[width = \linewidth]{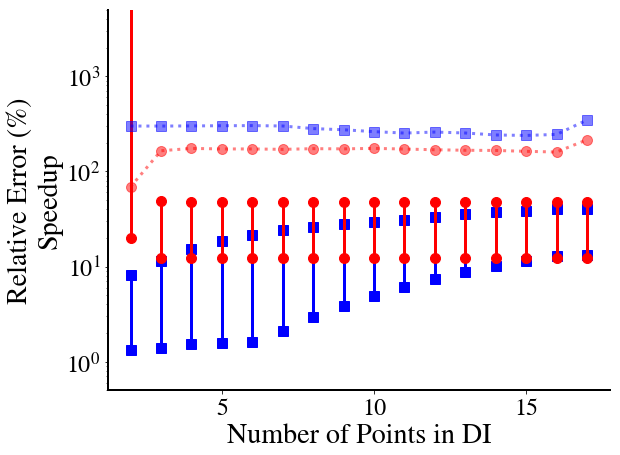}
        \caption{$\trainingSet = \{0.60,0.65,...,1.40\}$\\ $\testingSet = \{0.60,0.61,...,1.40\}$}
    \end{subfigure}
    \begin{subfigure}[b]{.45\linewidth}
        \centering
        \captionsetup{justification=centering}
        \includegraphics[width = \linewidth]{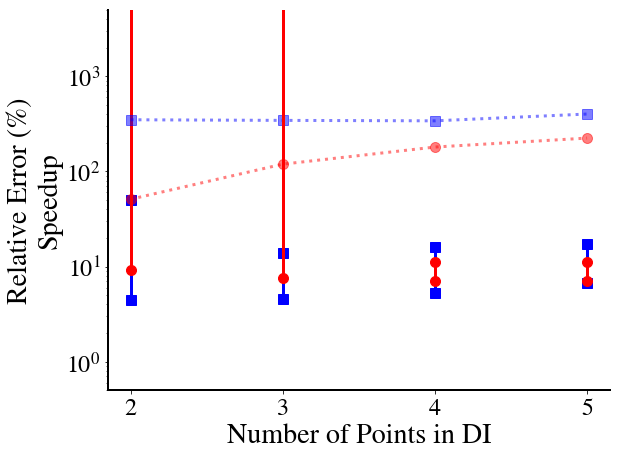}
        \caption{$\trainingSet = \{0.60,0.70,...,1.00\}$\\ $\testingSet = \{0.60,0.61,...,1.00\}$}
    \end{subfigure}
    \begin{subfigure}[b]{.45\linewidth}
        \centering
        \captionsetup{justification=centering}
        \includegraphics[width = \linewidth]{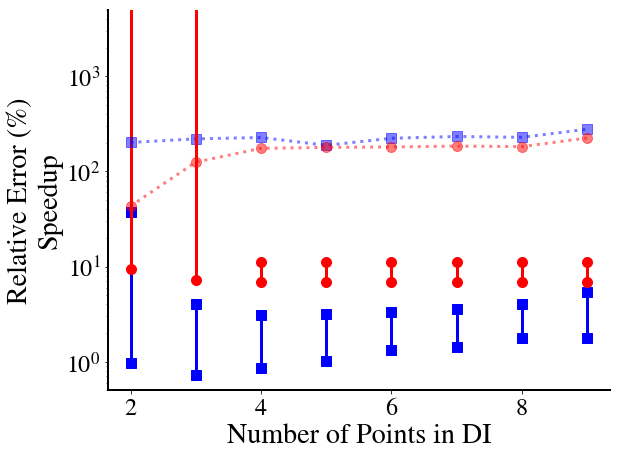}
        \caption{$\trainingSet = \{0.60,0.65,...,1.00\}$\\ $\testingSet = \{0.60,0.61,...,1.00\}$}
    \end{subfigure}
    \caption{The relative error ranges and speedups for the radial advection problem for all four sets of $\trainingSet$ described in Table~\ref{tab:training}. LaSDI-LS and LaSDI-NM use five and four latent-space dimensions, respectively. For $\trainingSet = \{0.60,0.65,\dots,1.0\}$, a quadratic dynamical system without cross-terms is used in DI. All other cases use the linear dynamical systems. For simplicity, we shorten the $y$-axis of each figure. The high error ranges correspond to numerically unstable dynamical systems found in DI.}
    \label{fig:ex9 error range}
\end{figure}

We now consider two different parameter spaces; a smaller one ($w\in[0.6,1.0]$) and a larger one ($w\in[0.6,1.4]$) as indicated in Table~\ref{tab:training}. We seek to analyze how flexible LaSDI can be in larger parameter spaces with potentially nonlinear latent space dynamics. Note that for this example, we only vary one parameter in the initial conditions. Thus, we do not include heat maps but rather graphs of the maximum relative error against the parameter value. We also vary the number of training points in each of these examples. Figure~\ref{fig:ex9 error range} shows the comparison of global and local DI that are trained on $\trainingSet$ and predicted on $\testingSet$. 

\begin{figure}
    \centering
    \begin{subfigure}[b]{.7\linewidth}
        \centering
        \includegraphics[width = \linewidth]{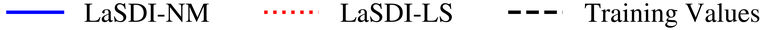}
    \end{subfigure}
    \begin{subfigure}[b]{.45\linewidth}
        \centering
        \includegraphics[width = \linewidth]{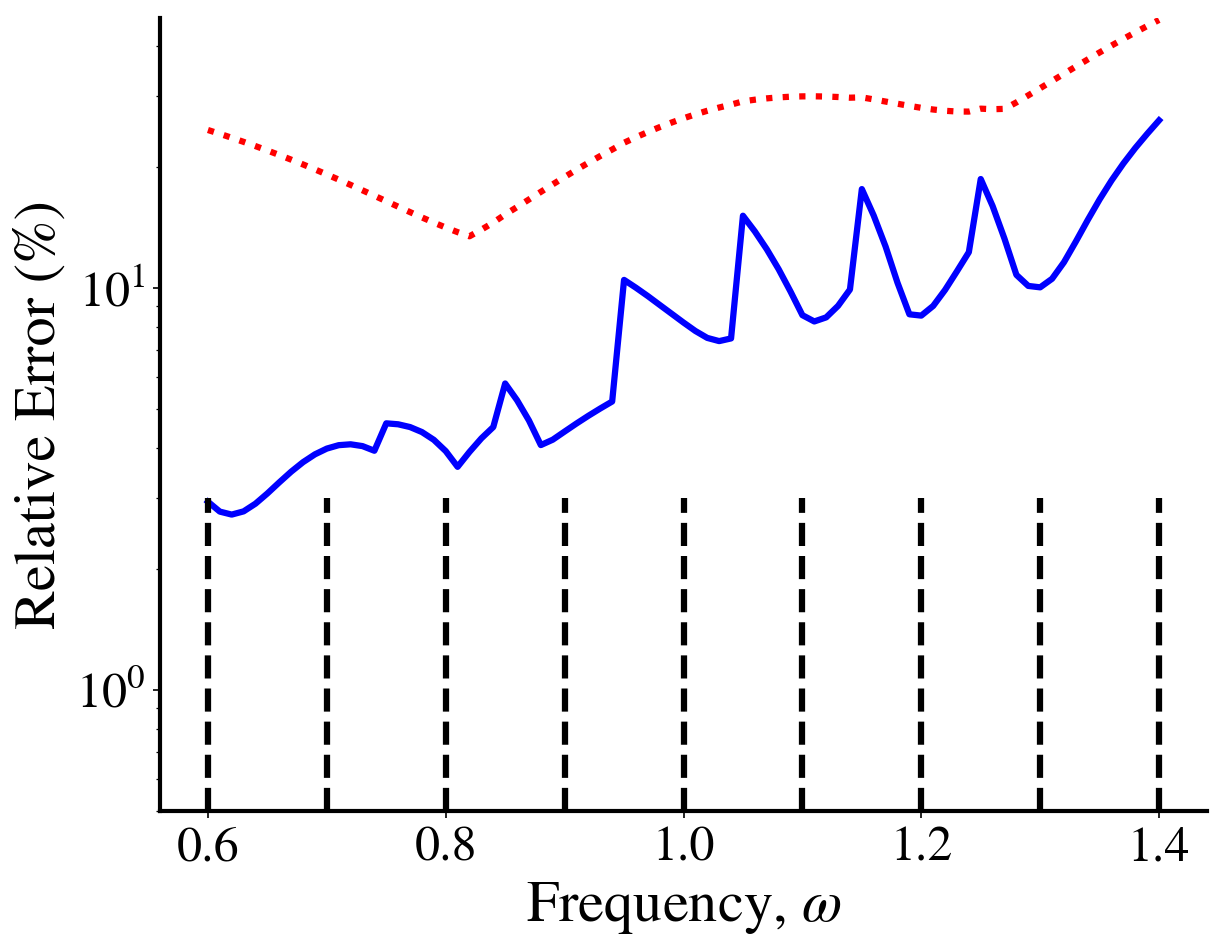}
        \captionsetup{justification=centering}
        \caption{Local (3) DI}
    \end{subfigure}
        \begin{subfigure}[b]{.45\linewidth}
        \centering
        \captionsetup{justification=centering}
        \includegraphics[width = \linewidth]{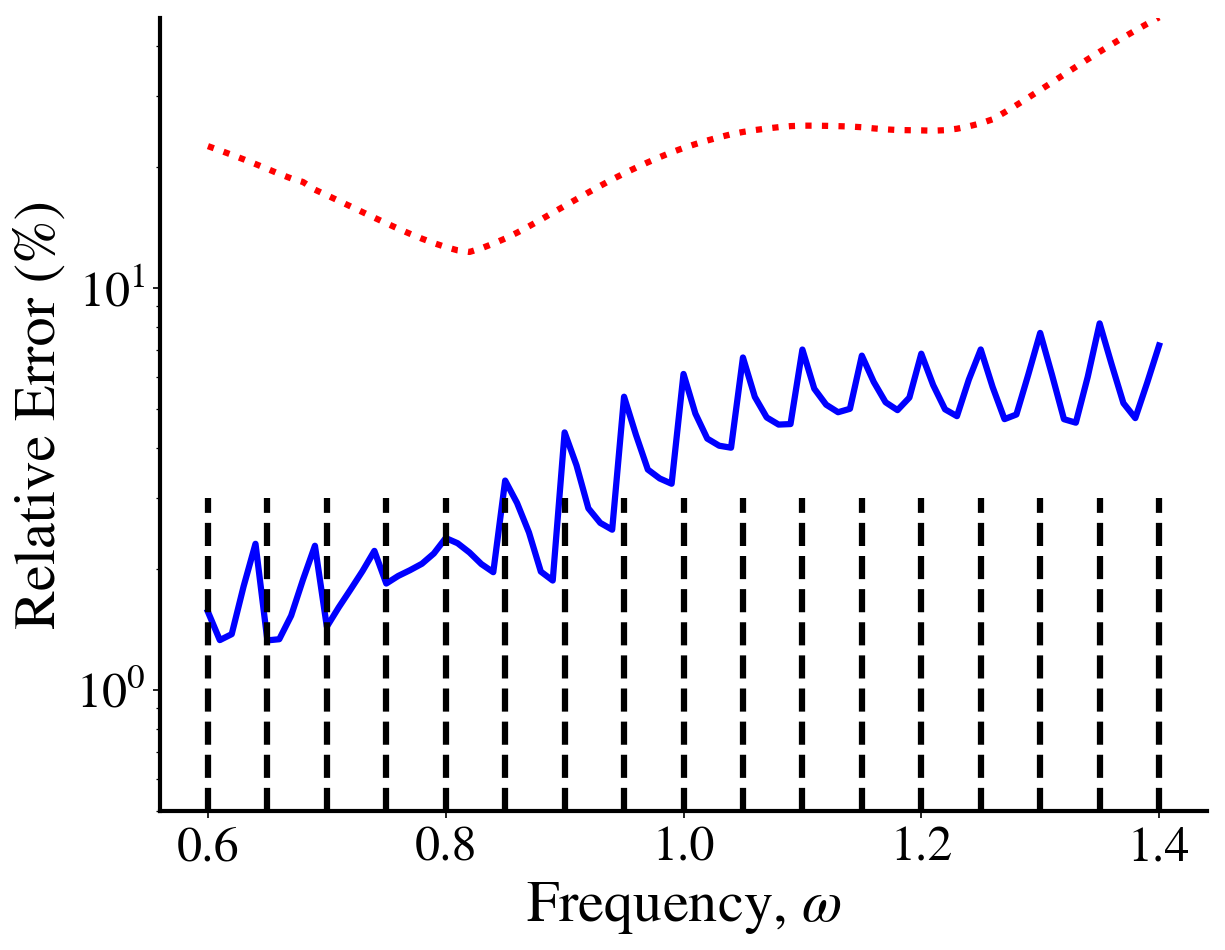}
        \caption{Local (2) DI}
    \end{subfigure}
    \begin{subfigure}[b]{.45\linewidth}
        \centering
        \captionsetup{justification=centering}
        \includegraphics[width = \linewidth]{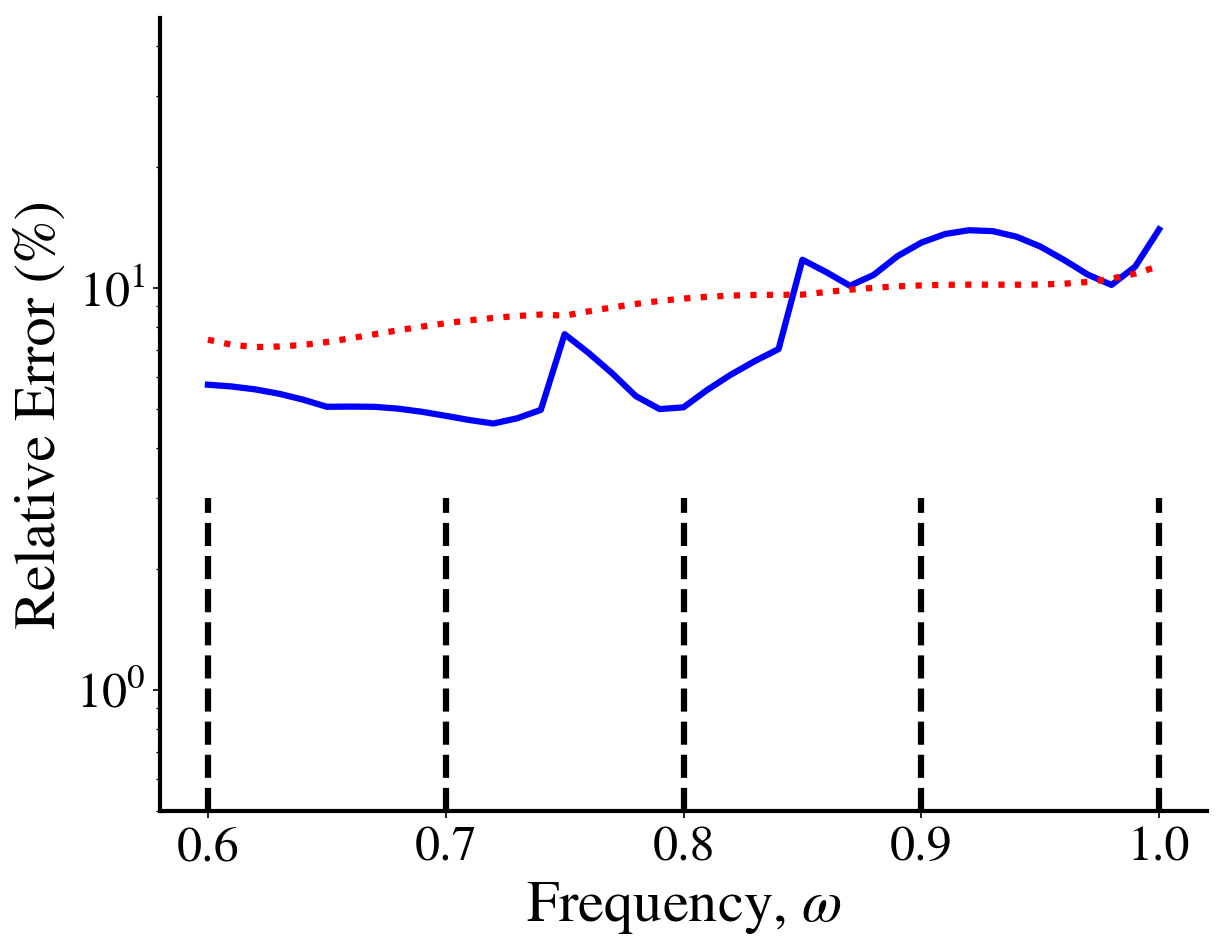}
        \caption{Local (3) DI}
    \end{subfigure}
    \begin{subfigure}[b]{.45\linewidth}
        \centering
        \captionsetup{justification=centering}
        \includegraphics[width = \linewidth]{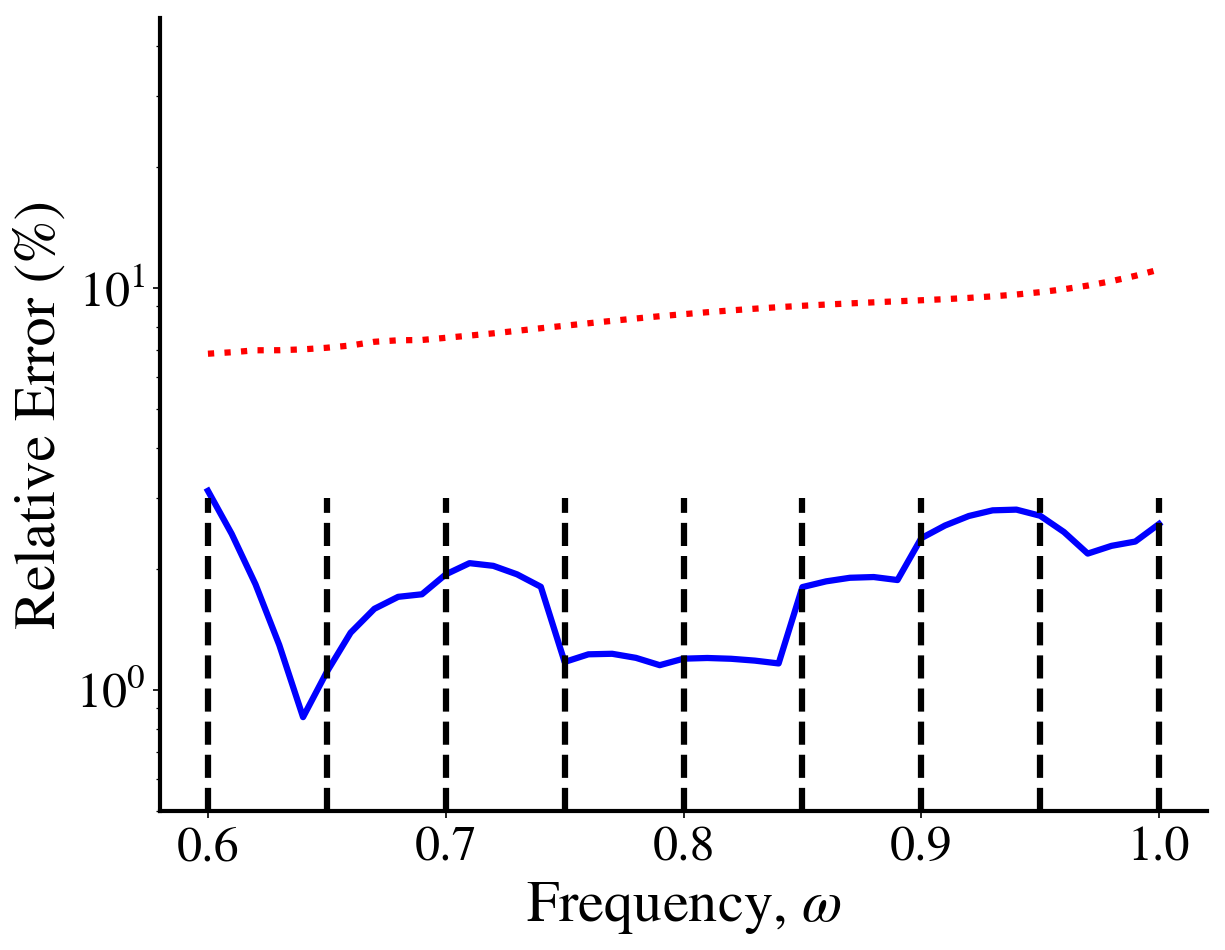}
        \caption{Local (4) DI}
    \end{subfigure}
    \caption{The relative error for the radial advection problem across the respective $\testingSet$ for both LaSDI-NM and LaSDI-LS with various local DI. Each example shows the best-case scenario from the error ranges illustrated in Figure~\ref{fig:ex9 error range}.}
    \label{fig:ex9 full}
\end{figure}

Figure \ref{fig:ex9 full} depicts the best-case scenario for each of these cases in terms of accuracy. We note the most accurate regime in each case is local DI. Among the best local LaSDI models, as expected, the smaller parameter space with the larger number of training points performs the best, with a maximum relative error of $\approx 3\%$. Likewise, the worst result appears in the larger parameter space with the smaller number of training points. In this case LaSDI-NM reaches error approximately 20\%.

Interestingly, the local LaSDI (2) DI model on $\trainingSet = \{0.60,0.65,\dots, 1.40\}$ (i.e., Figure~\ref{fig:ex9 full}(b)) outperforms the local LaSDI (3) DI model on $\trainingSet= \{0.60,0.70,\dots,1.00\}$ (i.e, Figure~\ref{fig:ex9 full}(c)) when using LaSDI-NM. This implies that local LaSDI models will perform well as long as there are enough training points whether the parameter space is large or small. On the other hand, the performance of the global LaSDI-NM models is affected by the parameter space size. For example, the worst error of the global LaSDI-NM for $\trainingSet = \{0.60,0.65,\dots, 1.40\}$ is $40.46\%$, while the worst error of the global LaSDI-NM for $\trainingSet= \{0.60,0.70,\dots,1.00\}$ is $17.18\%$.

\section{Conclusion}\label{sec:conclusion}
We have introduced a framework for latent space dynamics identification (LaSDI) that builds efficient parametric models by exploiting both local and global regression techniques. Two types of local dynamics identification models are introduced: (i) region-based and (ii) point-wise interpolated DIs. Latent spaces are generated in two different ways as well, i.e., using linear and nonlinear compression techniques. The compression techniques, such as proper orthogonal decomposition and auto-encoder, reduce large-scale simulation data to small-scale data whose size is manageable for regression techniques. The framework is general enough to be applicable to any physical simulations. We have applied various LaSDI models to four different problems, i.e., 1D and 2D Burgers, heat conduction, and radial advection problems. 

We observe that local versus global DI models need to be determined case by case. For example, a local DI model outperforms global DI for LaSDI-NM for 1D Burgers problem with nine training points (see the case of Degree 0 dynamical system in Figure~\ref{fig: 1db compare}) and 2D Burgers problem with 25 training points (see Figure~\ref{fig:2d GL compare}) and all the cases for radial advection problems (see Figure~\ref{fig:ex9 error range}), while a global DI or local DI with larger number of points outperform a local DI with smaller number of points for LaSDI-NM for 1D Burgers problem with 25 training points (see Figure~\ref{fig:1db sum}). 

The nonlinear compression techniques, such as auto-encoder, outperforms the linear compression technique for problems with slowly decaying Kolmogorov's width, while the linear compression techniques, such as proper orthogonal decomposition, outperforms the nonlinear compression technique for problems with fast decaying Kolmogorov's width. The nonlinear compression was able to bring the projection error down with small latent space dimension even for the problems with slowly decaying Komogorov's width. However, the computational cost of the nonlinear manifold compression is much higher than the linear compression, so we recommend using linear compression if the problem shows fast singular value decay. 

There are several future directions to consider in the context of latent space dynamics identification framework. First, we have only considered predefined uniformly distributed training parameters, which is fine with a small size parameter space dimension of one or two. However, the number of uniformly distributed training points increases exponentially as the dimension of the parameter space increases, which makes the generation of the simulation data extremely challenging. Therefore, adaptive and sparse parameter sampling for LaSDI needs to be developed, which will lead to an optimal sampling. Second, we have only considered the parameters that affects an initial condition. However, a parameter of interest may not alter the initial condition, e.g., material properties. The LaSDI framework introduced in this paper cannot handle this case. Therefore, parameterization other than the initial condition needs to be developed. Third, we have extensively used a regression technique to identify dynamics in a latent space. However, this is not the only option. For example, any system identification techniques or operator learning algorithms are applicable. Fourth, the architecture of the neural network for LaSDI-NM may not have been optimally chosen. The more thorough neural architecture search would lead to the better performance of LaSDI-NM models. Finally, because LaSDI is a non-intrusive method, it can be applied to data sets that might not be generated by simulations. These methods, along with adding regularization to the training process can help discover more robust latent-space dynamical systems for potentially noisy and imperfect data as done in \cite{Chen2021}.

With these future directions, the potential for LaSDI is to offer fast and accurate data-driven simulation capability. We expect to apply LaSDI to various important applications, such as climate science, manufacturing, and fusion energy. 

\section*{Acknowledgments}
This work was performed at Lawrence Livermore National Laboratory and partially funded by two LDRDs (21-FS-042 and 21-SI-006).  Lawrence Livermore National Laboratory is operated by Lawrence Livermore National Security, LLC, for the U.S. Department of Energy, National Nuclear Security Administration under Contract DE-AC52-07NA27344 and LLNL-JRNL-831849.

This research was supported in part by an appointment with the National Science Foundation (NSF) Mathematical Sciences
Graduate Internship (MSGI) Program sponsored by the NSF Division of Mathematical Sciences. This program is administered by the Oak Ridge Institute for Science and Education (ORISE) through an interagency agreement between the U.S. Department of Energy (DOE) and NSF. ORISE is managed for DOE by ORAU. All opinions expressed in this paper are the author's and do not necessarily reflect the policies and views of NSF, ORAU/ORISE, or DOE.

\bibliographystyle{ieeetr}
\bibliography{LaSDI_Bib.bib}
\end{document}